\documentclass[11pt]{amsart}

\def\C{\mathbb{C}}

\def\p{\varphi^{-1}}
\def\v{\varphi}
\def\l{\lambda}
\def\bl{\boldsymbol}
\def\o{\omega}
\def\O{\Omega}
\def\K{\mathcal K}
\def\D{\mathbb D}

\def\ov{\overline}
\def\lo{\longrightarrow}
\def\t{\tilde}

\def\mb{\mathbb}
\def\mr{\mathrm}
\def\w{\widehat}
\def\a{\alpha}
\def\b{\beta}
\def\g{\gamma}
\def\d{\displaystyle\sum}

\def\e{equivalent}
\def\eq{equivalence}
\def\mo{multiplication operator~}
\def\mos{multiplication operators~}
\def\i{\prime}

\def\u{unitarily equivalent}
\def\h{ Hermitian holomorphic vector bundle}
\def\ec{equivalence class}
\def\ue{unitary equivalence}
\def\rk{reproducing kernel~}
\def\rks{reproducing kernels~}

\def\c{corresponding}

\def\H{homogeneous}

\def\C{curvature~}
\def\Cs{curvatures~}

 \newtheorem{thm}{Theorem}[section]
 \newtheorem{cor}[thm]{Corollary}
 \newtheorem{lem}[thm]{Lemma}
 \newtheorem{prop}[thm]{Proposition}
 
 \newtheorem{rem}[thm]{Remark}
 \newtheorem{ex}[thm]{Example}
\newcommand{\be}{\begin{equation}}
\newcommand{\ee}{\end{equation}}
\newcommand{\bea}{\begin{eqnarray}}

\newcommand{\eea}{\end{eqnarray}}
\newcommand{\Bea}{\begin{eqnarray*}}
\newcommand{\Eea}{\end{eqnarray*}}

\newcounter{cnt1}
\newcounter{cnt2}
\newcounter{cnt3}
\newcommand{\blr}{\begin{list}{$($\roman{cnt1}$)$}
 {\usecounter{cnt1} \setlength{\topsep}{0pt}
 \setlength{\itemsep}{0pt}}}
\newcommand{\bla}{\begin{list}{$($\alph{cnt2}$)$}
 {\usecounter{cnt2} \setlength{\topsep}{0pt}
 \setlength{\itemsep}{0pt}}}
\newcommand{\bln}{\begin{list}{$($\arabic{cnt3}$)$}
 {\usecounter{cnt3} \setlength{\topsep}{0pt}
 \setlength{\itemsep}{0pt}}}
\newcommand{\el}{\end{list}}

\title[Curvature calculations for Homogeneous operators]{Curvature calculations for a class of \\homogeneous operators}

\author[Misra]{Gadadhar Misra}
\author[Shyam Roy]{Subrata Shyam Roy}
\address{Indian Statistical Institute\\
R. V. College Post\\
Bangalore 560 059\\
} \email[Gadadhar Misra]{gm@isibang.ac.in} \email[Subrata Shyam
Roy]{ssroy@isibang.ac.in}

\thanks{The research of the first author was supported in part by a
grant from the DST - NSF  Science and Technology Cooperation Programme.  
The second author was supported by the Indian Statistical Institute.}

\date{}
\begin{document}
\begin{abstract}
For an operator $T$ in the class ${\mathrm B}_n(\Omega)$, 
introduced in \cite{CD}, the simultaneous unitary equivalence
class of the curvature  and the covariant derivatives up to a
certain order of the corresponding bundle $E_T$  determine the 
unitary equivalence class of the operator $T$.  In the paper
\cite{CD2}, the authors ask if there exists some pair of 
inequivalent oprators $T_1$ and $T_2$ for which the 
simultaneous unitary equivalence class of the curvature 
along with all the covariant derivatives coincide except 
for the derivative of the highest order. Here we show that 
some of the covariant derivatives are necessary to determine 
the unitary equivalence class of the operators in 
${\mathrm B}_n(\Omega)$.  Our examples consist of
homogeneous operators. For homogeneous operators, the 
simultaneous unitary equivalence class of the curvature and 
all its covariant derivatives are determined from
the simultaneous unitary equivalence class of these at $0$. 
This shows that it is enough to calculate all the 
invariants and compare them at just one point, say $0$.
These calculations are then carried out in number of 
examples.
\end{abstract}

\maketitle
\section{Introduction}
For an open connected subset $\Omega$ of $\mathbb C$ and a
positive integer $n$, the class ${\mathrm B}_n(\Omega)$,
introduced in \cite{CD}, consists of bounded operators $T$ with the
following properties
\begin{enumerate}
\item[a)] $\Omega \subset \sigma (T)$ \item[b)] $\mbox{ran} (T-\omega)=
\mathcal H$ for $\omega \in\Omega$ \item[c)] $\bigvee_{\omega \in \Omega}
\ker(T-\omega)= \mathcal H$ for $\omega \in \Omega$ \item[d)] $\dim~\ker (T- \omega)=n$ for $\omega
\in\Omega$.
\end{enumerate}
A complete set of unitary invariants for the operators in the
class $\mr B_n(\O)$ was obtained in \cite{CD} as well. It was
shown in \cite[proposition 1.11]{CD} that the
 eigenspaces for each $T$ in  $\mathrm B_n (\Omega)$ form a
 Hermitian holomorphic vector bundle $E_T$ over $\Omega$, that is,
$$
E_T := \{(\o , x)\in \Omega \times \mathcal H : x \in \ker (T-
\o)\} \mbox{~and ~}\pi(\o , x)= \o
$$
and there exists a holomorphic frame $z\mapsto \gamma(\o)
:=(\gamma_1(\o), \ldots ,\gamma_n(\o) )$ with $\gamma_i(\o) \in
\ker (T-\o)$,  $1\leq i \leq n$.  The hermitian structure at $z$
is the one that $\ker (T-\o)$ inherits as a subspace of the
Hilbert space $\mathcal H$. In other words, the metric at $\o$ is
simply the grammian $h(\o) = \big ( \!\! \big (\langle
\gamma_j(\o), \gamma_i(\o) \rangle \big )\!\! \big )_{i,j=1}^n$.
The curvature $\mathcal K_T(\o)$ of the bundle $E_T$ is then
defined to be $\frac{\partial}{\partial \bar \o} \big ( h^{-1}
\frac{\partial}{\partial \o}h\big )(\o)$ for $\o\in \O$ (cf.
\cite[pp. 211]{CD}).
\begin{thm}[\cite{CD1}, Page. 326]
Two operators $T,\,\widetilde{T}$ in $\mathrm B_1(\Omega)$ are
unitarily equivalent if and only if $\mathcal K_T(\o)=\mathcal
K_{\widetilde{T}}(\o)$ for $\o$ in $\Omega$.
\end{thm}
Thus the  curvature of the line bundle $E_T$ is a complete set of
unitary invariant for an operator in $\mathrm B_1(\Omega)$.
It is not hard to see (cf. \cite[pp. 211]{CD}) that the curvature of a
bundle $E$ transforms according to the rule
$$\mathcal K(fg)(\o)=(g^{-1}\mathcal K(f)g)(\o), ~~~\o\in\Delta,$$
where $f=(e_1,...,e_n)$ is a frame for $E$ over an open subset
$\Delta \subseteq \Omega$ and $g:\Delta\longrightarrow
GL(n,\mathbb C)$ is a holomorphic map, that is, $g$ a holomorphic
change of frame. Since $g$ is a scalar valued holomorphic function
for a line bundle $E$, it follows from the transformation rule for
the curvature that it is independent of the choice of a frame in
this case. In general, the \C~ of a bundle $E$ of rank $n > 1$
depends on the choice of a frame. Thus the \C~ $\K(\o)$ itself
cannot be an invariant for the bundle $E$. However, the
eigenvalues of $\K(\o)$ are invariants for the bundle $E$. More
interesting is the description of a complete set of invariants
given in \cite{CD} involving the curvature and the covariant
derivatives $\K_{z^i{\bar z}^j}$, $0\leq i\leq j\leq i+j\leq n,
(i,j)\neq (0,n)$, where rank of $E=n$. They showed, in a
subsequent paper (cf. \cite{CD2}), by means of examples that fewer
covariant derivatives will not suffice to determine the class of
the bundle $E$. The examples they constructed do not necessarily
correspond to an operator of the class in $\mr B_n(\O)$. In this
paper we construct examples of operators $T$ in $\mr B_2(\D)$ and
$\mr B_3(\D)$ to show that the eigenvalues of \C~ alone does not
determine the class of the bundle $E_T$. Our examples show, one
will need at least derivatives of order $(1,1)$. Our examples
consists of bundles \H~ on the open unit disc $\D$. We will say
that a holomorphic Hermitian bundle $E$ over the unit disc $\D$ is
homogeneous if every bi-holomorphic automorphism $\varphi$ of the
unit disc lifts to an isometric isomorphism of the bundle $E$.
These verifications are somewhat nontrivial and use the
homogeneity of the bundle in an essential way. It is not clear if
for \H~ bundle the \C~ along with its derivatives up to order
$(1,1)$ suffices to determine its \ec. Secondly the original
question of sharpness of \cite[Page. 214]{CD} and \cite[page.
39]{CD2}, remains open, although our examples a provides partial
answer.

Let $B(z,w)=(1-z\bar w)^{-2}$ be the Bergman kernel on the unit
disc, the Hilbert space \c~ to the non-negative definite kernel
$B^{\l/2}(z,w)=(1-z\bar w)^{-\l}$ be $\mathbb A^{(\l)}(\D)$ for
$\l>0$. We let $M^{(\l)}:\mathbb A^{(\l)}(\D)\lo\mathbb
A^{(\l)}(\D)$ be the multiplication operator, that is,
$(M^{(\l)}f)(z)=zf(z),$ $f\in\mathbb A^{(\l)}(\D),z\in\D$.
Following the jet construction of \cite{DMV}, we construct a
Hilbert space $\mathbf A_k^{(\a,\b)}(\D)$ for $\a,\b>0$ consisting
of holomorphic functions defined on the open unit disc $\D$ taking
values in $\mathbb C^{k+1}$ starting from the kernel
$B^{(\a,\b)}(\mathbf z,\mathbf
w)=B^{\a/2}(z_1,w_1)B^{\b/2}(z_2,w_2)$, $\mathbf
z=(z_1,z_2),\mathbf w=(w_1,w_2)\in\D^2$. It turns
that the \rk for $\mathbf A_k^{(\a,\b)}(\D)$ is
$$B_k^{(\a,\b)}(z,w)={\big(\!\! \big( B^{\a/2}(z_1,w_1)\partial_{z_2}
\partial_{\bar w_2}B^{\b/2}(z_2,w_2)\big)\!\!\big)_{0\leq i,j\leq k}}_{|{\rm res}~ \D},$$
that is, $z_1=z=z_2$ and $w_1=w=w_2$.
The multiplictaion operator on $\mathbf A_k^{(\a,\b)}(\D)$ is
denoted by $M_k^{(\a,\b)}$.

For a suitably restricted class of operators, some times, the \ue~ class of the \C~ $\K_T$
determines the \ue~ class of the operator $T$. For instance, the
curvature at $0$ of the   generalised Wilkins operators
$M_k^{(\a,\b)}$ is of the form ${\rm diag}~\{\a, \cdots, \a, \a +
(k+1)\b + k(k+1)\}$. Thus the unitary equivalence class of the
curvature at $0$ determines the unitary equivalence class of these
operators within the class of 
the generalised Wilkins operarotrs of rank $k+1$ (cf. \cite{MSR}, \cite[page 428]{BM}).

\section{Examples from the Jet Construction}
\begin{ex} \rm
Consider the operators $M:=M^{(\l/2)}\oplus M^{(\mu/2)}$ and $
M^\i:=M_1^{(\a,\b)}$ for $\l,\mu>o$ and $\a,\b>0$. Wilkins
\cite{W} has shown that the operator $\tilde M^*$ is  in $\mr
B_2(\D)$ and that it is irreducible. This operator is also {\em
\H}, that is, $\v(\tilde M)$ is \ue~ to $\tilde M$ for all
bi-holomorphic automorphisms  $\v$ of the open unit disc $\D$ (cf.
\cite{BM}). It is easy to see that the operators $M^{(\lambda/2)}$
and $M^{(\mu/2)}$ are both homogeneous and the adjoint of these
operators are in the class $\mr B_1(\D)$. Consequently, the direct
sum, namely, $M^*$ is homogeneous and lies in the class $\mr
B_2(\D)$.  Let
\begin{enumerate}
\item $h(z)= \left ( \begin{array}{cc} B^{\l/2}(z,z) & 0 \\ 0 &
B^{\mu/2}(z,z) \end{array}\right )$, $\lambda,\,\mu >0$,
\item ${h}^\i(z)=B^{(\a,\b)}_1(z,z)^t=\left(%
\begin{array}{cc}
  (1-|z|^2)^2 & \b \bar z(1-|z|^2) \\
  \b z(1-|z|^2) & \b(1+\b |z|^2) \\
\end{array}%
\right)(1-|z|^2)^{-\a-\b-2}$, $\a,\b>0$, for $z\in  \mathbb D$.
\end{enumerate}
We see that $h$ and ${h}^\i$ are the metrics for bundles
corresponding to the the operators $M^*$ and ${M^\i}^*$
respectively. To emphasize the dependence of the \C~ on the
metric, we will find it useful to also write
$\K_h:=\bar\partial(h^{-1}\partial h)$.
\begin{enumerate}
\item[(a)] $\K_{h}(z)=$ diag
$(\l(1-|z|^2)^{-2},\mu(1-|z|^2)^{-2})$, \item[(b)]
$\big(\K_{h}\big)_{\bar z}(z) = 2~ {\rm diag} (\l
z(1-|z|^2)^{-3},\mu z(1-|z|^2)^{-3})$,
\item[(c)] $\K_{{h}^\i}(z)=\left(%
\begin{array}{cc}
  {\a} & {-2\b(\b+1)(1-|z|^2)^{-1}\bar z} \\
  {0} & {\a+2\b+2} \\
\end{array}%
\right)(1-|z|^2)^{-2}$,
\item[(d)] $\big(\K_{{h}^\i}\big)_{\bar z}(z)=2\left(%
\begin{array}{cc}
  {\a z} & {-\beta(\beta+1)(1+2|z|^2)(1-|z|^2)^{-1} } \\
  {0}& {(\a+2\b+2)z} \\
\end{array}%
\right)(1-|z|^2)^{-3} $.
\end{enumerate}
Choosing $\l>0$ and $\mu-\l>2$, we set $\a=\l$ and
$\b=\frac{1}{2}(\mu-\l-2)$. Since curvature is self-adjoint the
set of eigenvalues is the complete set of unitary invariants for
the curvature. The eigenvalues for $\mathcal K_{h}(z)$ and
$\mathcal K_{{h}^\i}(z)$, $z\in \mathbb D$, are clearly the same
by the choice of $\l,\mu,\a,\b$.    So these matrices are
pointwise unitarily equivalent.  Now we observe that
$\big(\K_{h}\big)_{\bar z}(0)=0$ and $\big(\K_{{h}^\i}\big)_{\bar
z}(0)\neq 0$. Hence they cannot be \u. Hence \C alone does not
determine the unitary \ec~ of the bundle.
\end{ex}
Before we construct the next example, let us recall that for any
\rk $K$ on $\D$, the normalized kernel $\tilde K(z,w)$ (in the
sense of Curto-Salinas \cite[Def.]{CS}) is defined to be the
kernel $K(0,0)^{1/2}K(z,0)^{-1}K(z,w)K(0,w)^{-1}K(0,0)^{1/2}$.
This kernel is characterized  by the property $\tilde K(z,0)=I$
and is therefore uniquely determined up to a conjugation by an
unitary matrix. Let $K(z,w)=\sum_{k,\ell\geq 0}a_{k\ell}z^k{\bar
w}^\ell$ and $\tilde K(z,w)=\sum_{k,\ell\geq 0}\tilde
a_{k\ell}z^k{\bar w}^\ell,$ where $a_{k\ell}$ and $\tilde
a_{k\ell}$ are determined by the real analytic functions $K$ and
$\tilde K$ respectively, $a_{k\ell}$ and $\tilde a_{k\ell}$ are in
$\mathcal M(n,\mb C)$, for $k,\ell\geq 0$. Since $\tilde K(z,w)$
is a normalized kernel, it follows that $\tilde a_{00}=I$ and
$\tilde a_{k0}=\tilde a_{0\ell}=0$ for $k,\ell\geq 1.$ Let
$K(z,w)^{-1}=\sum_{k,\ell\geq 0}b_{k\ell}z^k{\bar w}^\ell,$ where
$b_{k\ell}$ is in $\mathcal M(n,\mb C)$ for $k,\ell\geq 0$.
Clearly, $K(z,w)^*=K(w,z)$ for any \rk $K$ and $z,w\in\D$.
Therefore, $a_{k\ell}^*=a_{\ell k}$, ${\tilde a}_{k\ell}^*={\tilde
a_{\ell k}}$ and ${b_{k\ell}}^*=b_{\ell k}$ for $k,\ell\geq 0$,
where $X^*$ denotes the conjugate transpose of the matrix $X$.

If we assume that the adjoint $M^*$ of the \mo $M$ on the Hilbert
space $(\mathcal H,K)$ is in $\mr B_k(\D)$ then it is not  hard to
see that the operators $M^*$ on the Hilbert space $ \t {\mathcal
H}$ determined by the normalized kernel $\tilde K$ is \e~ to $M^*$
on the Hilbert space $(\mathcal  H,K)$. Hence the adjoint of the
\mo $M$ on $( \t{\mathcal  H}, \tilde K)$ lies in $\mr B_k(\D)$ as
well. Let $(E,\t h)$ be the \c~ bundle, where $\t h(z)=\tilde
K(z,z)^t$, $z\in \D$.
\begin{lem} \label{der}
If $\t h(z) = \tilde K(z,z)^t$, then ${\partial}^m \t
h(0)={\bar\partial}^n \t h(0)={\partial}^m \t
h^{-1}(0)={\bar\partial}^n \t h^{-1}(0)=0$ for $m,n\geq 1$ and
$\bar\partial\partial \t h(0)={\tilde a_{11}}^t$,
$\bar\partial\partial \t h^{-1}(0)=-{\tilde a_{11}}^t$,
${\bar\partial}^2{\partial}^2 \t h(0)=4{\tilde a_{22}}^t$.
\end{lem}
\begin{proof}
Let $\t h(z)=\d_{m,n\geq 0}h_{mn}z^m{\bar z}^n$ and $\t
h^{-1}(z)=\d_{m,n\geq 0}h^\prime_{mn}z^m{\bar z}^n$. As $\t
h(z)=\tilde K(z,z)^t$, $h_{mn}={\tilde a_{mn}}^t$ for $m,n\geq 0$,
so $h_{00}=I$, $h_{m0}=h_{0n}=0$ for $m,n\geq 1$. As
$h_{mn}=\frac{{\partial}^m{\bar\partial}^n\t h(0)}{m!n!}$ and
$h^\prime_{mn}=\frac{{\partial}^m{\bar\partial}^n{\t
h}^{-1}(0)}{m!n!}$ for $m,n\geq 0$, to prove the first assertion
it is enough to show that $h^\prime_{m0}=h^\prime_{0n}=0$ for
$m,n\geq 1$. As ${h^\prime_{k\ell}}^*=h^\prime_{\ell k}$ it is
enough to show that $h^\prime_{m0}=0$ for $m\geq 1$. It follows
from $\t h(z)\t h^{-1}(z)=I$ that $h^\prime_{00}=I$ and
$\d_{k=0}^{m}h_{m-k,0}h^\prime_{k0}=0$ for $m\geq 1$. So,
$0=\d_{k=0}^{m}h_{m-k,0}h^\prime_{k0}=\d_{k=0}^{m-1}h_{m-k,0}h^\prime_{k0}
+h^\prime_{m,0}=h^\prime_{m,0}$ for $m\geq 1$, as $h_{\ell 0}=0$
for $\ell\geq 1$.

For the second assertion we note that  $\t h(z)\t h^{-1}(z)=I$
implies $\bar\partial\partial \t
h^{-1}(0)=h^\prime_{11}=-h_{11}=-{\tilde a_{11}}^t$. Clearly
$\bar\partial\partial \t h(0)=h_{11}={\tilde a_{11}}^t$ and
${\bar\partial}^2{\partial}^2 \t h(0)=4h_{22}=4{\tilde a_{22}}^t$.
\end{proof}
\begin{lem} \label{cur}
The curvature $\K_{\t h}$ and the covariant derivative $(\K_{\t
h})_{{\bar z}^n}$of \C~ of the  bundle $(E,\t h)$ at $0$ is
$a_{11}^t$ and $(n+1)!a_{1,n+1}^t$, that is, $\K_{\t
h}(0)=a_{11}^t$ and $(\K_{\t h})_{{\bar
z}^n}(0)=(n+1)!a_{1,n+1}^t$.
\end{lem}
\begin{proof}
By \cite[page. 211]{CD} $\K_{\t h}=\bar\partial(\t h^{-1}\partial
\t h)$. Hence $\K_{\t h}(0)=\bar\partial \t h^{-1}(0)\partial \t
h(0)+\t h^{-1}(0)\bar\partial\partial \t h(0)={\tilde a}_{11}^t.$

For the second assertion, we know from \cite[Proposition 2.17,
page 211]{CD}, $\K_{{\bar
z}^n}(0)={\bar\partial}^n\K(0)={\bar\partial}^{n+1}(\t
h^{-1}\partial \t h )(0)$.  Using Leibnitz's rule, this is the
same as $\d_{k=0}^{n+1}\binom{n+1}{k} {\bar\partial}^{n+1-k}\t
h^{-1}(0){\bar\partial}^k\partial \t
h(0)={\bar\partial}^{n+1}\partial \t h(0)=
{(n+1)!}h_{1,n+1}=(n+1)!a_{1,n+1}^t$, as ${\bar\partial}^\ell \t
h^{-1}(0)=0$ for $\ell\geq 1$ from Lemma \ref{der}.
\end{proof}
\begin{lem} \label{dcur}
If $\K_{\t h}$ is the \C~ of the bundle $(E,\t h)$ then $(\K_{\t
h})_{z\bar z}(0)=2(2{\tilde a_{22}}-{{\tilde a_{11}}}^2)^t$.
\end{lem}
\begin{proof}
We know from \cite{CD} that for a bundle map of a \h~
$\Theta:(E,\t h)\lo(E,\t h)$ the covariant derivatives $\Theta_z$
and $\Theta_{\bar z}$ of $\Theta$ with respect to holomorphic
frame $f$ are $\Theta_ z(f)=\partial\Theta(f)+[\t h^{-1}\partial
\t h,\Theta(f)]$ and $\Theta_{\bar z}(f)=\bar\partial\Theta(f)$.
So $(\K_{\t h})_{z\bar z}(z)=\bar\partial\big(\partial(\K_{\t
h})(z)+[\t h^{-1}\partial \t h,\K_{\t
h}](z)\big)=\bar\partial\partial\K_{\t h}(z)+[\bar\partial(\t
h^{-1}\partial \t h),\K_{\t h}](z)+[\t h^{-1}\partial \t
h,\bar\partial\K_{\t h}](z)=\bar\partial\partial\K_{\t h}(z)+[\t
h^{-1}\partial \t h,\bar\partial\K_{\t h}](z).$ Hence $(\K_{\t
h})_{z\bar z}(0)=\bar\partial\partial\K(0),$ as $\partial \t
h(0)=h_{10}=0$. Now by Leibnitz rule we see that
$\bar\partial\partial\K_{\t
h}(z)=\bar\partial\partial(\bar\partial(\t h^{-1}\partial \t
h))(z)={\bar\partial}^2\big(\partial (\t h^{-1}\partial \t
h)\big)(z)={\bar\partial}^2\big(\partial \t h^{-1}\partial \t h+\t
h^{-1}\partial^2 \t h\big)(z)=\big({\bar\partial}^2\partial \t
h^{-1}(z)\partial \t h(z)+2\bar\partial\partial \t
h^{-1}(z)\bar\partial\partial \t h(z)+\partial \t
h^{-1}(z){\bar\partial}^2\partial \t
h(z)\big)+\big({\bar\partial}^2 \t h^{-1}(z)\partial^2 \t
h(z)+2\bar\partial \t h^{-1}(z)\bar\partial{\partial}^2 \t h(z)+\t
h^{-1}(z){\bar\partial}^2
\partial^2 \t h(z)\big)$. Therefore by Lemma \ref{der} we have
$(\K_{\t h})_{z\bar z}(0)=2\bar\partial\partial \t
h^{-1}(0)\bar\partial\partial \t h(0)+{\bar\partial}^2
\partial^2 \t h(0)=-2{\tilde a_{11}^t}{\tilde a_{11}^t}+4a_{22}^t
=2(2{\tilde a_{22}}-{{\tilde a_{11}}}^2)^t$.
\end{proof}
\begin{lem}\label{coeff}
The coefficient of $z^{k+1}{\bar w}^{\ell+1}$ in the power series
expansion of $\tilde K(z,w)$ is 
\begin{multline*}
\qquad  {\tilde
a}_{k+1,\ell+1} = a_{00}^{1/2}\big(\d_{s=1}^k\d_{t=1}^\ell b_{s0}
a_{k+1-s,\ell+1-t}b_{0t}+ \\   \d_{s=1}^kb_{s0}a_{k+1-s,\ell+1}b_{00}
+\d_{t=1}^\ell
b_{00}a_{k+1,\ell+1-t}b_{0t}+b_{00}a_{k+1,\ell+1}b_{00}-b_{k+1,0}a_{00}b_{0,\ell+1}\big)a_{00}^{1/2}
\end{multline*}
for $ k,\ell \geq 0$.
\end{lem}
\begin{proof}
From the definition of ${\tilde K}(z,w)$ we see that for
$k,\ell\geq 0$\Bea \lefteqn{ {\tilde a_{k+1,\ell+1} =
a_{00}^{1/2}\big(\d_{s=0}^{k+1}}\d_{t=0}^{\ell+1} b_{s0}
a_{k+1-s,\ell+1-t} b_{0t}\big)a_{00}^{1/2}}\\
&&=a_{00}^{1/2}\big(\d_{s=1}^{k+1}\d_{t=1}^{\ell+1} b_{s0}
a_{k+1-s,\ell+1-t} b_{0t} +\d_{s=1}^{k+1} b_{s0}
a_{k+1-s,\ell}b_{00}\\
&&\hspace*{72pt}+\d_{t=1}^{\ell+1} b_{00}a_{k+1,\ell+1-t} b_{0t}+
b_{00}a_{k+1,\ell+1}b_{00}\big)a_{00}^{1/2}\\
&&=a_{00}^{1/2}\big(\d_{s=1}^{k}\d_{t=1}^{\ell} b_{s0}
a_{k+1-s,\ell+1-t}
b_{0t}+\d_{s=1}^{k+1}b_{s0}a_{k+1-s,0}b_{0,\ell+1}+\d_{t=1}^\ell
b_{k+1,0}a_{0,\ell+1-t}b_{0t}\\
&&\hspace*{67pt}+\d_{s=1}^{k+1} b_{s0} a_{k+1-s,\ell}b_{00}
+\d_{t=1}^{\ell+1} b_{00}a_{k+1,\ell+1-t} b_{0t}+
b_{00}a_{k+1,\ell+1}b_{00}\big)a_{00}^{1/2}\\
&&=a_{00}^{1/2}\big(\d_{s=1}^{k}\d_{t=1}^{\ell} b_{s0}
a_{k+1-s,\ell+1-t}
b_{0t}+(\d_{s=0}^{k+1}b_{s0}a_{k+1-s,0})b_{0,\ell+1}+b_{k+1,0}(\d_{t=0}^{\ell+1}
a_{0,\ell+1-t}b_{0t})\\
&&\hspace*{30pt}+\d_{s=1}^{k} b_{s0} a_{k+1-s,\ell}b_{00}
+\d_{t=1}^{\ell} b_{00}a_{k+1,\ell+1-t} b_{0t}+
b_{00}a_{k+1,\ell+1}b_{00}-b_{k+1,0}a_{00}b_{0,\ell+1}\big)a_{00}^{1/2}\\
&&=a_{00}^{1/2}\big(\d_{s=1}^{k}\d_{t=1}^{\ell} b_{s0}
a_{k+1-s,\ell+1-t} b_{0t}+ \d_{s=1}^{k} b_{s0}
a_{k+1-s,\ell}b_{00} +\d_{t=1}^{\ell} b_{00}a_{k+1,\ell+1-t}
b_{0t}\\
&&\hspace*{83pt}+
b_{00}a_{k+1,\ell+1}b_{00}-b_{k+1,0}a_{00}b_{0,\ell+1}\big)a_{00}^{1/2}\\
\Eea as the coefficient of $z^{k+1} $ in $ K(z,w)^{-1}
K(z,w)=\sum_{s=0}^{k+1} b_{s0} a_{k+1-s,0}=0$ and the  coefficient
of ${\bar w}^{\ell+1} $ in $ K(z,w)
K(z,w)^{-1}=\sum_{t=0}^{\ell+1} a_{0,\ell+1-t} b_{0t}=0$ for
$k,\ell \geq 0.$
\end{proof}
The following Theorem will be useful in the sequel.
For $T$ in $\mr B_n(\O)$,  recall that $\K_T$ denotes the \C~ of the bundle
$E_T$ \c~  to $T$.
\begin{thm} \label{hom}
Suppose that $T_1$ and $T_2$ are \H~ operators in $\mathrm
B_n(\D)$ then $\K_{T_1}(0)$ and $(\K_{T_1})_{\bar z}(0)$ are
simultaneously \u~ to $\K_{T_2}(0)$ and $(\K_{T_2})_{\bar z}(0)$
respectively if and only if $\K_{T_1}(z)$ and $(\K_{T_1})_{\bar
z}(z)$ are simultaneously \u~ to $\K_{T_2}(z)$ and
$(\K_{T_1})_{\bar z}(z)$ respectively for $z$ in $\D$.
\end{thm}
Before going into the proof of \ref{hom} let us fix some
notations. Let M\"{o}b denote the group of biholomorphic
automorphisms on the unit disc $\D$ in the complex plane,
$c:~~$M\"{o}b$\times\mathbb D\longrightarrow\mathbb C$ be the
function which is given by the formula
$c(\varphi^{-1},z):=(\varphi^{-1})^\prime(z)$, where the prime
stands for differentiation with respect to $z$. The function $c$
satisfies the following {\it cocycle} property:
$$ c(\varphi^{-1}\psi^{-1},z)=c(\varphi^{-1},\psi^{-1}(z))c(\psi^{-1},z), \mbox{~for~}
\varphi\in \mbox{M\"{o}b and}~ z\in \mathbb D.$$ This cocyle
property can easily verified by chain rule.
\begin{lem} \label{hcurv}
Suppose that $T$ in $\mathrm B_n(\mathbb D)$ is \H.  Then \Bea
&&(\mr{a})~ \mathcal
K_T(\varphi^{-1}(0))=|c(\varphi^{-1},0)|^{-2}U^{-1}_\varphi
\mathcal K_T(0)U_\varphi \\ &&  (\mr b)~(\mathcal K_T)_{\bar
z}(\varphi^{-1}(0))= |c(\varphi^{-1},0)|^{-2}
\ov{{c(\varphi^{-1},0)}^{-1}}U^{-1}_\varphi \Big((\mathcal
K_T)_{\bar z}(0)-
\ov{c(\varphi^{-1},0)^{-1}(\varphi^{-1})^{(2)}(0)}\mathcal
K_T(0)\Big)U_\varphi \Eea

for some unitary
operator $U_\varphi$, $\varphi \in$ {\rm M\"{o}b}.

\end{lem}
\begin{proof}   From \cite{CD} it follows that homogeneity of $T$ implies $\mathcal
K_{\varphi(T)}(z)=U^{-1}_{\varphi,z}\mathcal K_T(z)U_{\varphi,z}$
for some unitary operator $U_{\varphi,z}$, $\varphi \in \mbox{\rm
M\"{o}b}$ and $z \in \mathbb D$. Note that an application of chain
rule gives the formula

\begin{equation} \label{chainrule}
\K_{\varphi(T)}(z) = |(\varphi^{-1})^\prime(z)|^2\mathcal
K_T((\varphi^{-1})(z)), \mbox{~for~}\v \in \mbox{\rm M\"{o}b}
\mbox{~and~} z \in \mathbb D.
\end{equation}
Now, assuming that $\mathcal K_{\varphi(T)}(z)=
U^{-1}_{\varphi,z}\mathcal K_T(z)U_{\varphi,z}$, using
\eqref{chainrule}, amounts to $U^{-1}_{\varphi,z} \mathcal
K_T(z)U_{\varphi,z}
= |c(\varphi^{-1},z)|^{2} \mathcal K_T((\varphi^{-1})(z)).$
Putting $z=0$, we get $U^{-1}_{\varphi,0}\mathcal
K_T({\varphi}^{-1}(0))U_{\varphi,0} =
|c(\varphi^{-1},0)|^{2}\mathcal K_T((\varphi^{-1})(0))$. If we set
$U_{\varphi,0}:=U_\varphi$, then we have $\mathcal
K_T({\varphi}^{-1}(0))=|c(\varphi^{-1},0)|^{-2}U^{-1}_\varphi\mathcal
K_T(0)U_\varphi$~~for $\varphi \in$ {\rm M\"{o}b}, $z\in \mathbb
D$. This proves (a).

To prove (b), we recall \eqref{chainrule} and differentiating with
respect to $z$ we get
\begin{eqnarray}
\bar\partial\mathcal K_{\varphi(T)}(z)&=&
{(\varphi^{-1})^\prime(z)}\overline{(\varphi^{-1})^{(2)}(z)}
\mathcal K_T(\varphi^{-1}(z))+|(\varphi^{-1})^\prime(z)|^2
\overline{(\varphi^{-1})^\prime(z)}\bar\partial\mathcal
K_T(\varphi^{-1}(z))\nonumber \\ \label{transrule}
&=&{c(\varphi^{-1},z)}\overline{(\varphi^{-1})^{(2)}(z)}\mathcal
K_T(\varphi^{-1}(z))+|c(\varphi^{-1},z)|^2
\overline{c(\varphi^{-1},z)}\bar\partial\mathcal
K_T(\varphi^{-1}(z)).
\end{eqnarray}
Using (\ref{transrule}) and (a), putting $z=0$  and
$U_{\varphi,0}=U_\varphi$, we see that
\begin{eqnarray}
\lefteqn{U^{-1}_\varphi\bar\partial\mathcal K_T(0)U_\varphi =
{c(\varphi^{-1},0)}\overline{(\varphi^{-1})^{(2)}(0)}\mathcal
K_T(\varphi^{-1}(0))+|c(\varphi^{-1},0)|^2
\overline{c(\varphi^{-1},0)}\bar\partial\mathcal K_T(\varphi^{-1}(0))} \nonumber \\
&&={c(\varphi^{-1},0)}\overline{(\varphi^{-1})^{(2)}(0)}|c(\varphi^{-1},0)|^{-2}
U^{-1}_\varphi\mathcal K_T(0)U_\varphi+|c(\varphi^{-1},0)|^2
\overline{c(\varphi^{-1},0)}\bar\partial\mathcal
K_T(\varphi^{-1}(0)) \nonumber\\ \label{trans0}
&&=\overline{c(\varphi^{-1},0)^{-1}(\varphi^{-1})^{(2)}(0)}U^{-1}_\varphi\mathcal
K_T(0)U_\varphi+|c(\varphi^{-1},0)|^2
\overline{c(\varphi^{-1},0)}\bar\partial\mathcal
K_T(\varphi^{-1}(0)).
\end{eqnarray}
So, $\bar\partial\mathcal K_T(\varphi^{-1}(0))=
|c(\varphi^{-1},0)|^{-2}\ov{c(\varphi^{-1},0)^{-1}}U^{-1}_\varphi
\Big(\bar\partial\mathcal
K_T(0)-\ov{c(\varphi^{-1},0)^{-1}(\varphi^{-1})^{(2)}(0)}\mathcal
K_T(0)\Big)U_\varphi$. As $(\K_T)_{\bar z}=\bar\partial\K_T$ by
\cite{CD}, this proves (b).
\end{proof}
\begin{cor} \label{one}
Suppose that $T_1$, $T_2$ are \H~ operators in $\mathrm B_n(\D)$.
Then \begin{enumerate} \item[(1)]$U^{-1}\K_{T_2}(0)U=\K_{T_1}(0)$
\item[(2)] $U^{-1}(\K_{T_2})_{\bar z}(0)U=(\K_{T_1})_{\bar z}(0)$
\end{enumerate} for some
unitary operator $U$ if and only if
\begin{enumerate}
\item[(i)]$V_{\v}^{-1}\K_{T_2}(z)V_{\v}=\K_{T_1}(z)$
\item[(ii)]$V_{\v}^{-1}(\K_{T_2})_{\bar
z}(z)V_{\v}=(\K_{T_1})_{\bar z}(z)$
\end{enumerate} for some unitary operator  $V_{\v}$,
$\v\in$M\"{o}b, $z\in\D$.
\end{cor}
\begin{proof}
One part is obvious, let us prove the other part.

Take $\v=\v_{t,z}$, where $\v_{t,a}(z)=t\frac{z-a}{1-\bar az}$,
for $a,z,z\in\D$ and $t\in\mb T$. Pick a unitary operator such
that $(a)$ and (b) of Lemma \ref{hcurv} are satisfied. We get from
(1) and Lemma \ref{hcurv}(a) that \Bea
\K_{T_2}(z)&=&|c(\p,0)|^{-2}U^{-1}_{\varphi} \K_{T_2}(0)U_\varphi \\
 &=& |c(\p,0)|^{-2}U^{-1}_{\varphi} U^{-1}\K_{T_1}(0)U
U_\varphi  \\
&=& |c(\p,0)|^{-2}U^{-1}_{\varphi}
U^{-1}|c(\p,0)|^{2}U_\v\K_{T_1}(z)U^{-1}_\v U U_\varphi \\
\label{a} &=&U^{-1}_{\varphi} U^{-1}U_\v\K_{T_1}(z)U^{-1}_\v U
U_\varphi. \Eea As the product of unitary operators is again a
unitary operator, taking $V_\v=U^{-1}_\v U U_\varphi$ we have (i).

From $(2)$ and Lemma \ref{hcurv}(b)  \bea
\lefteqn{(\K_{T_2})_{\bar z} (z)=
|c(\p,0)|^{-2}\ov{c(\p,0)^{-1}}U^{-1}_\varphi \Big((
\K_{T_2})_{\bar z}(0)-\ov{c(\p,0)^{-1}(\p)^{(2)}(0)} \K_{T_2}(0)\Big)U_\varphi}\nonumber \\
 &=&|c(\p,0)|^{-2}\ov{c(\p,0)^{-1}}U^{-1}_\v U^{-1}
\Big(( \K_{T_1})_{\bar
z}(0)-\ov{c(\p,0)^{-1}(\p)^{(2)}(0)}\mathcal
K_{T_1}(0)\Big)U U_\v\nonumber\\
&=& |c(\p,0)|^{-2}\ov{c(\p,0)^{-1}}U^{-1}_\v
U^{-1} \Big(\ov{(\p)^{(2)}(0)c(\p,0)^{-1}}\K_{T_1}(0)\nonumber\\
 &&+|c(\p,0)|^{2}\ov{c(\p,0)}U_\v( \K_{T_1})_{\bar z}(z)U^{-1}_\v
-\ov{(\p)^{(2)}(0)c(\p,0)^{-1}}\K_{T_1}(0)\Big)U U_\v\nonumber\\
\label{b} &=&U^{-1}_\v U^{-1} U_\v(\K_{T_1})_{\bar z}(z)U^{-1}_\v
U U_\v. \eea  Taking $V_\v=U^{-1}_\v U U_\varphi$  as before, we
have (ii).
\end{proof}
\begin{proof}[Proof of Theorem \ref{hom}]
Combining Lemma \ref{hcurv} and Corollary \ref{one} we have a
proof of the Theorem \ref{hom}.
\end{proof}
For a positive integer $m$ let $S_m(c_1,\ldots,c_m)$ denote the
forward shift on $\mb C^{m+1}$ with weight sequence
$\{c_1,\ldots,c_m\}$, $c_i\in\mb C$, that is,
$\big(\!\!\big(S_m(c_1,\ldots,c_m)\big)\!\!\big)_{\ell,
p}=c_\ell\delta_{p+1,\ell}$, $0\leq p,\ell\leq m$. We set $\mb
S_m:=S_m(1,\ldots,m)$.
\begin{ex} \rm
Consider the operators $M_1=M^{(\a/2)}\oplus M_1^{(\a,\b^\prime)}$
and $M_2=M_2^{(\a,\b)}$ for $\a,\b,\b^\prime>0$. Wilkins \cite{W}
has shown that the adjoint of the operator $
{M_1^{(\a,\b^\prime)}}$ is in $\mr B_2(\D)$. This operator is also
\H. It is easy to see that the operator $M^{(\a/2)}$ is
homogeneous and its adjoint is in the class $\mr B_1(\D)$.
Consequently, the direct sum, namely, ${M_1}^*$ is homogeneous and
lies in the class $\mr B_3(\D)$. The operator ${M_2}^*$ is in $\mr
B_3(\D)$ by \cite[Proposition 3.6]{DMV} and is \H~ by \cite[Page.
428]{BM} and \cite[Theorem 5.1]{MSR}. Let
 $h_1(z)=(1-|z|^2)^{-\a}\oplus
B_1^{(\a,\b^\prime)}(z,z)^t$ and $h_2(z)=B_2^{(\a,\b)}(z,z)^t$. We
see that $h_1$ and $h_2$ the metrics for the bundles   $E_1$ and
$E_2$ \c~ to the operators $M_1^*$
and $M_2^*$ respectively, where $B_1^{(\a,\b^\prime)}(z,w)=\left(%
\begin{array}{cc}
  (1-z\bar w)^2 & \b^\prime  z(1-z\bar w) \\
  \b^\prime\bar w(1-z\bar w) & \b^\prime(1+\b^\prime z\bar w) \\
\end{array}%
\right)(1-z\bar w)^{-\a-\b^\prime-2}$ and $B_2^{(\a,\b)}(z,w)=\left(%
\begin{array}{ccc}
  (1-z\bar w)^4 & \beta(1-z\bar w)^3 z & \beta(\beta+1)(1-z\bar w)^2{ z}^2 \\
  \beta(1-z\bar w)^3 \bar w & \beta(1+\beta z\bar w)(1-z\bar w)^2 & \beta(\beta+1)(2+\beta z\bar w)(1-z\bar w) z \\
  \beta(\beta+1)(1-z\bar w)^2{\bar w}^2 & \beta(\beta+1)(2+\beta z\bar w)(1-z\bar w) 
\bar w&\beta(\beta+1)(2+(\beta+1)(4+ \beta z\bar w)z\bar w) \\
\end{array}%
\right)\\(1-z\bar w)^{-\alpha-\beta-4}$.
\begin{lem}\rm \label{relcur}
The \C~ at zero and the covariant derivatives of \C~ at zero up to
order for the bundles $E_1$ and $E_2$ respectively are
\begin{enumerate}
\item[(a)]$\K_1(0)=\mr{diag}(\a,\a,\a+2\b^\prime+2)$ ,
$(\K_1)_{\bar
z}(0)=S_2(0,-2\sqrt{\beta^\prime}\big(\beta^\prime+1)\big)^t $ and
$(\K_1)_{z\bar z}(0)= 2$ $\mr{diag}
(\a,\a+\b^\prime(\b^\prime+1),\a+\b^\prime(-\b^\prime+1)+2)$
\item[(b)]$\K_2(0)=\mr{diag}~(\a,\a,\a+3\b+6)$, $(\K_2)_{\bar
z}(0)=S_2\big(0,-{3}\sqrt{2(\beta+1)}(\beta+2)\big)^t$ and
$(\K_2)_{z\bar
z}(0)=\mr{diag}\big(\a,\a+3(\b+1)(\b+2),\a-3\b(\b+2)\big)$,
\end{enumerate} where $\K_i$, $(\K_i)_{\bar
z}$ and $(\K_i)_{z\bar z}$ are computed with respect to a metric
normalized at $0$ obtained from $h_i$ for $i=1,2$, that is, with
respect to an orthonormal basis at $0$.
\end{lem}
\begin{proof}
For any \rk $K$  with $K(z,w)=\d_{m,n\geq 0}a_{mn}z^m{\bar w}^n$
and $K(z,w)^{-1}=\d_{m,n\geq 0}b_{mn}z^m{\bar w}^n$the identity
$K(z,w)^{-1}K(z,w)=I$ implies that $b_{00}=a_{00}^{-1}$ and
$\d_{\ell=0}^kb_{0,k-\ell}a_{0\ell}=0$ for $k\geq 1$. For $k=1$ we
have $b_{10}=-a_{00}^{-1}a_{10}a_{00}^{-1}$, $b_{01}=(b_{10})^*$.
We have by Lemma \ref{coeff} \bea\lefteqn{\tilde
a_{11}=a_{00}^{1/2}\big(b_{00}a_{11}b_{00}-b_{10}a_{00}b_{01}\big)a_{00}^{1/2}}\nonumber\\\label{id1}
&&=a_{00}^{-1/2}\big(a_{11}-a_{10}a_{00}^{-1}a_{01}\big)a_{00}^{-1/2}.
\eea For $k=2$ we have
$b_{02}=-\big(b_{01}a_{01}+b_{00}a_{02}\big)a_{00}^{-1}=a_{00}^{-1}\big(a_{01}a_{00}^{-1}
a_{01}-a_{02}\big)a_{00}^{-1}$. We get from Lemma \ref{coeff}
\bea\lefteqn{\tilde
a_{12}=a_{00}^{1/2}\big(b_{00}a_{11}b_{01}+b_{00}a_{12}b_{00}
-b_{10}a_{00}b_{02}\big)a_{00}^{1/2}}\nonumber\\\label{id2}
&&=a_{00}^{-1/2}\big(a_{12}
-(a_{11}-a_{10}a_{00}^{-1}a_{01})a_{00}^{-1}a_{01}-a_{10}a_{00}^{-1}a_{02}\big)a_{00}^{-1/2}.\eea
Observing that
$b_{20}={b_{02}}^*=a_{00}^{-1}\big(a_{10}a_{00}^{-1}
a_{10}-a_{20}\big)a_{00}^{-1}$, from Lemma \ref{coeff} we have
\bea\lefteqn{\tilde a_{22}=
a_{00}^{1/2}\big(b_{10}a_{11}b_{01}+b_{10}a_{12}b_{00}+
b_{00}a_{21}b_{01}+b_{00}a_{22}b_{00}-b_{20}a_{00}b_{02}\big)a_{00}^{1/2}}\nonumber\\
&&=a_{00}^{-1/2}\big(a_{10}a_{00}^{-1}a_{11}a_{00}^{-1}a_{01}-a_{10}a_{00}^{-1}a_{12}
-a_{21}a_{00}^{-1}a_{01}+a_{22}\nonumber\\
&&\hspace*{83pt}-(a_{10}a_{00}^{-1}a_{10}-a_{20})a_{00}^{-1}
(a_{01}a_{00}^{-1}a_{01}-a_{02})\big)a_{00}^{-1/2}\nonumber\\\label{id3}
&&=a_{00}^{-1/2}\big(a_{22}+
(a_{20}a_{00}^{-1}a_{01}-a_{21})a_{00}^{-1}a_{01}-a_{20}a_{00}^{-1}a_{02}\nonumber
\\
&&\hspace*{80pt}-a_{10}a_{00}^{-1}(a_{12}
-(a_{11}-a_{10}a_{00}^{-1}a_{01})a_{00}^{-1}a_{01}
-a_{10}a_{00}^{-1}a_{02})\big)a_{00}^{-1/2}.\eea  We get from
$h_1$ that
$a_{11}=\mr{diag}~(\a,\a+\b^\prime,\b^\prime(\a+2\b^\prime+2))$,
$a_{00}=\mr{diag}(1,1,\b^\prime)$,
$a_{10}=S_2\big(0,\b^\prime\big)^t$,
$a_{12}=S_2\big(0,\b^\prime(\a+\b^\prime+1)\big) $,
$a_{22}=\mr{diag}\big(\frac{\a(\a+1)}{2},\frac{(\a+\b^\prime)(\a+\b^\prime+1)}{2}
,\frac{\b^\prime(\a+\b^\prime+2)(\a+3\b^\prime+3)}{2}\big)$,
$a_{02}=0$. So,
$a_{11}-a_{10}a_{00}^{-1}a_{01}=\mr{diag~(\a,\a,\b^\prime(\a+2\b^\prime+2))}$,
hence from Lemma \ref{cur} and Equation \eqref{id1} we have
$\K_1(0)={\tilde a_{11}}^t=\mr{diag}~(\a,\a,(\a+2\b^\prime+2))$.

We get from Equation \eqref{id2} $\tilde
a_{12}=S_2\big(0,-\sqrt{\b^\prime}(\b^\prime+1)\big)$, so from
Lemma \ref{cur} we have $(\K_1)_{\bar z}(0)=2\tilde
a_{12}^t=2S_2\big(0,-\sqrt{\b^\prime}(\b^\prime+1)\big)^t$.

From Equation \eqref{id3} $\tilde
a_{22}=\mr{diag}\big(\frac{\a(\a+1)}{2},\frac{\a(\a+1)+\b^\prime(\b^\prime+1)}{2},
\frac{(\a+\b^\prime+2)(\a+3\b^\prime+3)}{2}\big)$, hence from
Lemma \ref{dcur} we get $(\K_1)_{z\bar z}(0)=2(2{\tilde
a_{22}}-{{\tilde a_{11}}}^2)^t=$ 2
 $\mr{diag}(\a,\a+\b^\prime(\b^\prime+1),\a+\b^\prime(-\b^\prime+1)+2)$.
 This completes the proof of (a).

 To prove (b) we get from $h_2$ that
 $a_{00}=\mr{diag}(1,\b,2\b(\b+1))$, $a_{10}=S_2\big(\b, 2\b(\b+1)\big)^t$,
 $a_{12}=S_2\big(\b(\a+\b+1),\b(\b+1)(2\a+3\b+6)\big)$, $$\big(a_{02}\big)_{ij}=\left\{%
\begin{array}{ll}
    \b(\b+1), & {\mr{for}~i=3,j=1;} \\
    0, & \hbox{otherwise.} \\
\end{array}%
\right.$$ $a_{11}$ is a diagonal matrix with diagonal entries
$\a+\b,\b(\a+2\b+2)$ and $ 2\b(\b+1)(\a+3\b+6)$ respectively and
$a_{22}$ is also a diagonal matrix with diagonal entries
$\frac{(\a+\b)(\a+\b+1)}{2},\frac{\b(\a+\b+2)(\a+3\b+3)}{2}$
 and ${\b(\b+1)
\big((\a+\b+4)(\a+\b+5)+4(\b+1)(\a+\b+4)+\b(\b+1)\big)}$
respectively. Therefore
$a_{11}-a_{10}a_{00}^{-1}a_{01}=\mr{diag}(\a,\a\b,2\b(\b+1)(\a+3\b+6))$,
hence from Lemma \ref{cur} and Equation \eqref{id1} we have
$\K_2(0)={\tilde a_{11}}^t=\mr{diag}(\a,\a,\a+3\b+6)$.

We get from Equation \eqref{id2},  $\tilde
a_{12}=S_2\big(0,-\frac{3}{\sqrt 2}\sqrt{\b+1}(\b+2)\big)$. Now
from Lemma \ref{cur} we have $(\K_2)_{\bar z}(0)=2\tilde
a_{12}^t=S_2\big(0,-3\sqrt {2(\b+1)}(\b+2)\big)^t$. From Equation
\eqref{id3}, $\tilde
a_{22}=\mr{diag}\big(\frac{\a(\a+1)}{2},\frac{\a(\a+1)+3(\b+1)(\b+2)}{2}
,\frac{\a(\a+1)}{2}+3(\b+2)(\a+\b+3)\big)$, using Lemma \ref{dcur}
we get $(\K_2)_{z\bar z}(0)=2(2{\tilde a_{22}}-{{\tilde
a_{11}}}^2)^t=$ 2
$\mr{diag}\big(\a,\a+3(\b+1)(\b+2),\a-3\b(\b+2)\big)$.
\end{proof}
We prove a sequence of lemmas which  exhibits a unitary between
the vector spaces $\big((E_1)_0, h_1(0)\big)$ and $\big((E_2)_0,
h_2(0)\big)$ which intertwines $\K_1(0),\K_2(0)$ and $(\K_1)_{\bar
z}(0),(\K_2)_{\bar z}(0)$, where $(E_1)_0$ and $(E_2)_0$ are the
fibres over $0$ of the corresponding bundles.
\begin{lem}\label{uni}
  $U_0:(\mathbb C^3,h_2(0))\lo(\mathbb C^3,h_1(0)),$ is a
diagonal unitary with $U_0=\mbox{\rm diag\,}(\a_1,\a_2,\a_3),$ $\a_i\in\mathbb
C$ for $i=1,2,3,$ if and only if $|\a_1|^2=1, |\a_2|^2=\beta,
|\a_3|^2=\frac{2\beta(\beta+1)}{\beta^\prime}.$
\end{lem}
\begin{proof}
``only if'' part:  As $U_0$ is a unitary $U_0^*=U_0^{-1},$ where
$*$ denotes the adjoint of $U_0.$ Now, from \cite[p. 395]{DMV}  
\Bea
U_0^* &=& h_2(0)^{-1}\overline U_0^t h_1(0)\\
&=& \mr{diag}\big(1,\beta^{-1},(2\beta(\beta+1))^{-1}\big)
\mr{diag}(\bar \a_1,\bar \a_2,\bar \a_3)
\mr{diag}(1,1,\beta^\prime)\\
&=&\mr{diag}\big(\bar \a_1,\frac{\bar \a_2}{\beta},
\frac{\bar \a_3\beta^\prime}{2\beta(\beta+1)}\big)\\
&=& \mr{diag}(\a_1^{-1},\a_2^{-1},\a_3^{-1}) 
\Eea
This implies the desired equalities.

``if'' part: Taking $\alpha_1=1,\alpha_2=\sqrt\beta,\alpha_3=
 \sqrt{\frac{2\beta(\beta+1)}{\beta^\prime}},$ we see that $U_0=\mbox{diag}(\alpha_1,
 \alpha_2,\alpha_3)$ is a unitary between the two given vector
 spaces.
\end{proof}

The proof of the next lemma is just a routine verification.
\begin{lem}\label{inter}
Suppose that $T$ and $\tilde T$ are in $\mathcal M(3,\mb C)$ such that
$(T)_{ij}=\left\{%
\begin{array}{ll}
    \eta, & {{for}~i=2,j=3;} \\
    0, & \hbox{otherwise.} \\
\end{array}%
\right.$     and $\tilde T=\left\{%
\begin{array}{ll}
    \tilde\eta, & {{for}~i=2,j=3;} \\
    0, & \hbox{otherwise.} \\
\end{array}%
\right.    $ are two matrices, $T$ and $\tilde T$ satisfies
$AT=\tilde TA$ for some invertible diagonal matrix
$A=\mbox{diag}(a_1,a_2,a_3)$ if and only if
$\frac{\tilde\eta}{\eta}=\frac{a_2}{a_3}$.
\end{lem}
\begin{lem}\label{inter1}
If $\b^\prime=\frac{3}{2}\b+2$ then $U_0^{-1}\K_1(0)U_0=\K_2(0)$
and $U_0^{-1}(\K_1)_{\bar z}(0)U_0=(\K_2)_{\bar z}(0),$ where
$U_0:\big(\mathbb C^3,h_2(0)\big)\lo\big(\mathbb C^3,h_1(0)\big)$,
is a diagonal unitary with $U_0=diag(\alpha_1,\alpha_2,\alpha_3),$
$\alpha_i\in\mathbb C$ for $i=1,2,3$.
\end{lem}
\begin{proof}
By the choice of $\b^\prime$, $\K_1(0)=\K_2(0)$ by Lemma
\ref{relcur}, so  the first equality is clear.

 Take $T=(\K_2)_{\bar z}(0)$ and $\tilde T=(\K_1)_{\bar z}(0)$. Also choose
  $\alpha_1=1,\alpha_2=\sqrt\beta,\alpha_3=
 \sqrt{\frac{2\beta(\beta+1)}{\beta^\prime}}$, with $\b^\prime=\frac{3}{2}\b+2$.  To complete
the proof of the second equality, by Lemma \ref{inter} we only
have to verify $\frac{\tilde
\eta}{\eta}=\frac{\alpha_2}{\alpha_3}$, where
$\eta=-3\sqrt{2(\beta+1)}(\beta+2),\tilde\eta=-2\sqrt{\beta^\prime}(\beta^\prime+1)$
and $\alpha_1,\alpha_2,\alpha_3$ as above. Now
$\frac{\alpha_2}{\alpha_3}=\sqrt{\frac{\beta\beta^\prime}{2\beta(\beta+1)}}
=\sqrt{\frac{\frac{3}{2}\beta+2}{2(\beta+1)}}=\frac{1}{2}\sqrt{\frac{3\beta+4}{\beta+1}}$
and $\frac{\tilde
\eta}{\eta}=\frac{-2\sqrt{\beta^\prime}(\beta^\prime+1)}{-3\sqrt{2(\beta+1)}(\beta+2)}
=\frac{2{\sqrt{\frac{3}{2}\beta+2}{(\frac{3}{2}\beta+2+1)}}}{3\sqrt{2(\beta+1)}(\beta+2)}
=\frac{3\sqrt{3\beta+4}(\beta+2)}{2.3\sqrt{\beta+1}(\beta+2)}
=\frac{1}{2}\sqrt{\frac{3\beta+4}{\beta+1}}.$ Hence we have proved
the lemma.
\end{proof}
As the operators $M_1$ and $M_2$ are \H, Combing Lemma
\ref{inter1} with Theorem \ref{hom} we have the following
\begin{cor}\label{both}
There exists a unitary operator $U_\v$ such that
${U_\v}^{-1}\K_1(z)U_\v=\K_2(z)$ and ${U_\v}^{-1}(\K_1)_{\bar
z}(z)U_\v=(\K_2)_{\bar z}(z)$ for $z$ in $\D$, where $\v=\v_{t,z}$
in \mbox{\rm M\"{o}b} for $(t,z)\in \mb T\times \D$.
\end{cor}
\begin{lem}  \label{ineq}
If $\b^\prime=\frac{3}{2}\b+2$ then $(\K_1)_{z\bar z}(0)$ and
$(\K_2)_{z\bar z}(0)$ are not \u.
\end{lem}
\begin{proof}
By Lemma \ref{relcur} $(\K_i)_{z\bar z}(0)=\mr{diag}(p_i,q_i,r_i)$
for $i=1,2$, where $p_1=\a$, $q_1=\a+\b^\prime(\b^\prime+1)$,
$r_1=\a+\b^\prime(-\b^\prime+1)$ and $p_2=\a$,
$q_2=\a+3(\b+1)(\b+2)$, $r_2=\a-3\b(\b+2)$. As
$\b^\prime=\frac{3}{2}\b+2$, $q_1=\a+\frac{3}{4}(\b+2)(3\b+4)$ and
$r_1=\a-\frac{1}{4}(3\b+2)(3\b+4)$. So clearly $p_1=p_2$,
$q_1>r_1$ and  $q_2>r_2$  As $(\K_1)_{z\bar z}(0)$ and
$(\K_2)_{z\bar z}(0)$ are diagonal matrices, they are \u~ if and
only if $p_1=p_2$, $q_1=q_2$ and $r_1=r_2$. We see that $q_1\neq
q_2$ and  $r_1\neq r_2$, hence $(\K_1)_{z\bar z}(0)$ and
$(\K_2)_{z\bar z}(0)$ are not \u.
\end{proof}
Hence we have proved the following Theorem.
\begin{thm}
The simultaneous \ue~ class of the \Cs and the covariant
derivatives of the \Cs of order $(0,1)$ for the operators $M_1$
and $M_2$ are the same for $\beta^\prime = \tfrac{3}{2} \beta +2$. 
However, the covariant derivatives of the \Cs of order $(1,1)$ are 
not \u.
\end{thm}
\end{ex}
\section{Irreducible Examples and Permutation of Curvature Eigenvalues}
 In the first example constructed above one of the two \H~ operators
$M^*$ is reducible while the other ${\tilde M}^*$ is irreducible.
Similarly in the second example one of the two operators $M_1^*$
is reducible whereas the other $M_2^*$ is irreducible.
Irreducibility of ${\tilde M}^*$ and $M_2^*$ follows from
\cite{MSR}. We are interested in constructing such examples within
the class of irreducible operators in $\mathrm B_n(\D)$. The class
of irreducible homogeneous operators in $\mathrm B_2(\D)$ cannot
possibly possess such examples. Therefore, we consider a class of
\H~ operators in $\mr B_3(\D)$ discussed in \cite{KM}.

Let $\lambda$ be a real number and $m$ be a positive integer such
that $2\lambda-m>0.$ For brevity, we will write
$2\l_j=2\lambda-m+2j,$ $0\leq j\leq m$. Let
$$L(\lambda)_{\ell
j}=\left\{%
\begin{array}{ll}
    \binom{\ell}{j}^2
\frac{(\ell-j)!}{(2\l_j)_{\ell-j}}, & {\mr{for}~ 0\leq j\leq \ell\leq m;} \\
    0, & \hbox{otherwise.} \\
\end{array}%
\right.$$ and  $\mathrm B=$diag $(d_0,d_1,\ldots,d_m)$. Now
consider for $\bl\mu=(\mu_0,\ldots,\mu_m)$ with $\mu_0=1$ and
$\mu_\ell>0$ for $\ell=1,\ldots, m$
$$K^{(\lambda,\boldsymbol\mu)}(z,w)=(1-z\bar w)^{-2\lambda-m}D(z\bar w)\mbox{~exp}(\bar w
\mb S_m)\mathrm{B}\mbox{~exp}(z\mb S_m^*)D(z\bar w),$$ where
$\mathrm{B}_{\ell,\ell}=d_\ell=\d_{j=0}^\ell\binom{\ell}{j}^2
\frac{(\ell-j)!}{(2\l_j)_{\ell-j}}\mu_j^2$  for $0\leq \ell\leq
m$, that is, $L(\lambda)\boldsymbol\mu^\prime=\boldsymbol d$ for
$\bl\mu^\prime=(\mu_0^2,\mu_1^2,\ldots,\mu_m^2)^t$  and
$\boldsymbol d=(d_0,d_1,\ldots,d_m)^t$, $D(z\bar w)=(1-z\bar
w)^{m-\ell}\delta_{p\ell}$ is diagonal and $\mb S_m$ is the
forward shift with weight sequence $\{1,\ldots,m\},$ that is,
$\big(\!\!\big(\mb
S_m\big)\!\!\big)_{\ell,p}=\ell\delta_{p+1,\ell}$, $0\leq
p,\ell\leq m$, $X^t$ denotes the transpose of the matrix $X$.
$K^{(\lambda,\boldsymbol\mu)}$ is the reproducing kernel for the
Hilbert space $\mathbf A^{(\lambda,\bl\mu)}(\mathbb D)$ of
$\mathbb C^{m+1}$-valued holomorphic functions described in
\cite{KM}. Let $M^{(\l,\bl\mu)}$ denote the \mo on the Hilbert
space $\mathbf A^{(\lambda,\bl\mu)}(\mathbb D)$. In \cite{KM} it
is shown that $M^{(\l,\bl\mu)}$ is \H~ and irreducible, moreover,
${M^{(\l,\bl\mu)}}^*$ is in $\mr B_{m+1}(\D)$.
\begin{lem} \label{curv}
For the \rk $K^{(\l,\bl\mu)}$
\begin{enumerate}
\item[(a)]$\tilde a_{11}=[\mathrm B^{-1}\mb S_m\mathrm B,\mb
S_m^*]+(2\lambda+m)I_{m+1}-2D_m$, \item[(b)]$\tilde a_{12}={\mr
B}^{1/2}\big(\frac{1}{2}({\mr B}^{-1}\mb S_m^2\mr B{\mb S_m^*}{\mr
B}^{-1}+{\mb S_m^*}{\mr B}^{-1}{\mb S_m}^2)+{\mr B}^{-1}[D_m,\mb
S_m]-{\mr B}^{-1}\mb S_m\mr B\mb S_m^*{\mr B}^{-1}\mb S_m\big){\mr
B}^{1/2}$.
\end{enumerate}

where $I_k$ denotes the identity matrix of order $k$ and
$D_m=\mr{diag}~(m,\ldots,1,0).$
\end{lem}
\begin{proof}
From Equation \eqref{id1} in Lemma \ref{relcur} we get $\tilde
a_{11}=a_{00}^{-1/2}\big(a_{11}-a_{10}a_{00}^{-1}a_{01}\big)a_{00}^{-1/2}.$
Form the expansion of the reproducing kernel $K^{(\l,\bl\mu)}$ we
see that $a_{00}=\mathrm B$, $a_{10}=\mathrm B\mb S_m^*$,
$a_{01}=\mb S_m\mathrm B$, $a_{11}=\mb S_m\mathrm B\mb
S_m^*+(2\lambda+m)\mathrm B-2D_m\mathrm B$. So,
$a_{11}-a_{10}a_{00}^{-1}a_{01}=\mb S_m\mathrm B\mb
S_m^*+(2\lambda+m)\mathrm B-2D_m\mathrm B-\mathrm B\mb
S_m^*\mathrm B^{-1}\mb S_m\mathrm B.$ The proof of (a) is now
complete since the matrices  $\mb S_m\mathrm B\mb S_m^*$, $\mb
S_m{\mathrm B}^{-1}\mb S_m^*$, $\mathrm B$, $\mr B^{1/2}$, $\mr
B^{-1/2}$  are diagonal .

From Lemma \ref{coeff}, we have  $\tilde
a_{12}=a_{00}^{1/2}\big(b_{00}a_{11}b_{01}+b_{00}a_{12}b_{00}
-b_{10}a_{00}b_{02}\big)a_{00}^{1/2}$. Again, from the expansion
of the \rk $K^{(\l,\bl\mu)}$ it is easy to see that
$a_{12}=\frac{1}{2}\mb S_m^2\mr B {\mb S_m^*}+(2\l+m)\mb S_m\mr B
- D_m \mb S_m\mr B-\mb S_m\mr BD_m$, $b_{00}={\mr B}^{-1}$,
$b_{10}=-\mb S_m^*{\mr B}^{-1}$, $b_{02}=\frac{1}{2}{\mr
B}^{-1}{\mb S_m}^2$. To complete the proof of (b), it is enough to
note that two diagonal matrices $\mr B$ and $D_m$ commute.
\end{proof}

  It will be convenient to let  $\K_{\l,\bl\mu}$
 denote the \C~
 $\K_{\t h}(z)=\frac{\partial}{\partial\bar z}\big({\t h}^{-1}\frac{\partial }
 {\partial z}\t h\big)(z)$, where $\t h(z)=\tilde
 K^{(\l,\bl\mu)}(z,z)^t$ for $z$ in $\D$.
  Recall that $\tilde
 K^{(\l,\bl\mu)}$ is  the normalized \rk obtained from the \rk
 $K^{(\l,\bl\mu)}$.  Now we specialize to the case $m=2$.
\begin{lem}\label{curvp}
The \C~ at zero $\K_{\l,\bl\mu}(0)$ and the  covariant derivative
of \C~ at zero $(\K_{\l,\bl\mu})_{\bar z}(0)$ are
\begin{enumerate}
\item[(a)]$\K_{\l,\bl\mu}(0)=\mr{diag}~(a-b-2,a+b-c,a+c+2)$,
\item[(b)]$(\K_{\l,\bl\mu})_{\bar z}(0)=2S_2\big(-\sqrt
b(1+b-\frac{c}{2}),-\sqrt c(1+c-\frac{b}{2})\big)^t$, where
$a=2\lambda, b=d_1^{-1}, c=4d_1d_2^{-1}$.
\end{enumerate}

\end{lem}
\begin{proof}
$\mathrm B^{-1}\mb S_2\mathrm B\mb
S_2^*=\mbox{diag}(0,d_1^{-1},4d_1d_2^{-1})$ and $\mb S_2^*\mathrm
B^{-1}\mb S_2\mathrm B=\mbox{diag}(d_1^{-1},4d_1d_2^{-1},0)$.
Therefore by Lemma \ref{curv}(a) we see that $\tilde
a_{11}=\mbox{diag}~(2\lambda-d_1^{-1}-2,2\lambda+d_1^{-1}-4d_1d_2^{-1},2\lambda+4d_1d_2^{-1}+2)
=\mbox{diag}~(a-b-2,a+b-c,a+c+2).$ Hence $\K_{\l,\bl\mu}(0)=\tilde
a_{11}^t=\mbox{diag}~(a-b-2,a+b-c,a+c+2)$.

For (b) we note that ${\mr B}^{-1}\mb S_2^2\mr B{\mb S_2^*}{\mr
B}^{-1}=S_2(0,2d_1^{-1}d_2^{-1})$, $\mb S_2^*{\mr B}^{-1}{\mb
S_2}^2=S_2(4d_2^{-1},0)$, ${\mr B}^{-1}[D_2,\mb
S_2]=S_2(-d_1^{-1},-2d_2^{-1})$, ${\mr B}^{-1}\mb S_2\mr B\mb
S_2^*{\mr B}^{-1}\mb S_2=S_2( d_1^{-2}, 8d_1d_2^{-2})$ we have the
desired conclusion from Lemma \ref{cur} and Lemma \ref{curv}(b).
\end{proof}
If $\delta_1,\delta_2,\delta_3$ are the eigenvalues $\K(0)$ then
we know from \cite[Proposition 2.20]{CD} that $\delta_i>0$ for
$i=1,2,3$. Now, suppose $(\delta_1,\delta_2,\delta_3),$ is  a
fixed ordered triple of positive numbers. Then there exists
$K^{(\lambda,\boldsymbol \mu)}$ with $\l>1$ and $\mu_\ell>0$ for
$\ell=1,2$ such that
$\K_{\l,\bl\mu}(0)=\mbox{diag}(\delta_1,\delta_2,\delta_3)$,  only
if $\delta_i$'s satisfy the inequalities of Lemma \ref{nece}
below.
\begin{rem} \label{tuple}
We  emphasize that the \rk~ $K^{(\lambda,\boldsymbol\mu)}$ is
computed from a ordered basis, that is,
$K^{(\lambda,\boldsymbol\mu)}(w,w)=\big ( \!\! \big (\langle
\gamma_i(w), \gamma_j(w) \rangle \big )\!\! \big )_{i,j=1}^{3}$,
where $\{\g_1(w),\g_2(w),\g_3(w)\}$ is an ordered basis.
Consequently, the eigenvalues of $\K_{\l,\bl\mu}(0)$, which is
diagonal, appear in a fixed order. If one considers
$\{\g_{\sigma(1)}(w),\g_{\sigma(2)}(w),\g_{\sigma(3)}(w)\}$, it
will give rise to a different \rk~
$P_{\sigma}K^{(\lambda,\boldsymbol\mu)}{P_{\sigma}}^{*}$, say
$K^{(\lambda,\boldsymbol\mu)}_\sigma$, where $\sigma\in \Sigma_3$,
$\Sigma_3$ denotes the symmetric group of degree $3$ and
$$ (P_\sigma)_{i,j}=\left\{%
\begin{array}{ll}
    1, & {for~(i,j)=(i,\sigma(i));} \\
    0, & \hbox{otherwise.} \\
\end{array}%
\right.$$ Hence,
$\K_{h_\sigma}(0)=\mr{diag}(\delta_{\sigma(1)},\delta_{\sigma(2)},\delta_{\sigma(3)})$,
where ${h_\sigma}(z)=\tilde
K_\sigma^{(\lambda,{\boldsymbol\mu})}(z,z)^t$. It follows that the
curvature of the \c~ bundle as a matrix depends on the choice of
the particular ordered basis. The set of eigenvalues of \C~ at
$0$, which is diagonal in our case, will be thought of as an
ordered tuple, namely the ordered set of diagonal elements of
$\K_{\l,\bl\mu}(0)$.
\end{rem}
\begin{lem} \label{para}
$(\lambda,\bl\mu)=(\lambda^\prime,\bl\nu)$ if and only if
$(a,b,c)=(a^\prime,b^\prime,c^\prime)$ as ordered tuples,  where
$\bl\mu=(1,\mu_1,\mu_2)$, $\bl\nu=(1,\nu_1,\nu_2)$,
$\mu_\ell,\nu_\ell>0$ for $\ell=1,2$, $\bl d=(1,d_1,d_2)^t$, $\bl
d^\prime=(1,{d_1}^\prime,{d_2}^\prime)^t$; for $2\g_j=2\g-2+2j$,
$\g=\lambda,\lambda^\prime$$$L(\g)_{\ell
j}=\left\{%
\begin{array}{ll}
    \binom{\ell}{j}^2
\frac{(\ell-j)!}{(2\g_j)_{\ell-j}}, & {\mr{for}~ 0\leq j\leq \ell\leq m;} \\
    0, & \hbox{otherwise.} \\
\end{array}%
\right.$$  $\bl d=L(\lambda){\bl\mu}^\prime$, ${\bl
d}^\prime=L(\lambda^\prime){\bl\nu}^\prime$,
${\bl\mu}^\prime=(1,{\mu_1}^2,{\mu_2}^2)^t$,
${\bl\nu}^\prime=(1,{\nu_1}^2,{\nu_2}^2)^t$ $0\leq j\leq i\leq 2$,
$a=2\lambda$, $b={d_1}^{-1}$, $c=4d_1{d_2}^{-1}$,
$a^\prime=2\lambda^\prime$, $b^\prime={{d_1}^\prime}^{-1}$,
$c^\prime=4{d_1}^\prime{{d_2}^\prime}^{-1}$.
\end{lem}
\begin{proof}
One implication  is clear, so prove the other implication.
$a=a^\prime$ implies that $\lambda=\lambda^\prime$. $b=b^\prime$
and $c=c^\prime$ imply that $\bl d={\bl d}^\prime$. Now
invertibility of $L(\lambda)$ implies that
${\bl\mu}^\prime={\bl\nu}^\prime$, that is, $\bl\mu=\bl\nu$.
\end{proof}
\begin{lem} \label{nece}
If $(\delta_1,\delta_2,\delta_3)$ is an ordered tuple of positive
numbers such that $\K_{\l,\bl\mu}(0)=$ diag
$(\delta_1,\delta_2,\delta_3)$ then  $$\begin{array}{c}
  \delta_1+\delta_2+\delta_3>6 \\
   \delta_2+\delta_3-2\delta_1>6 \\
  2\delta_3-\delta_1-\delta_2>6.   \\
\end{array}$$
\end{lem}
\begin{proof}
By Lemma \ref{curvp}, Remark \ref{tuple} and the hypothesis of the
Lemma we have
$$\begin{array}{c}
a-b-2=\delta_1,\\
a+b-c=\delta_2,\\
a+c+2=\delta_3. \\
\end{array}$$

Equivalently, $A\bl x=\bl b,$ where $A=\left(%
\begin{smallmatrix}
  1 & -1 & 0 \\
  1 & 1 & -1 \\
  1 & 0 & 1 \\
\end{smallmatrix}%
\right),$ $\bl x=\left(%
\begin{smallmatrix}
  a \\
  b \\
  c \\
\end{smallmatrix}%
\right),$ $\bl {b} =\left(
\begin{smallmatrix}
  \delta_1+2 \\
  \delta_2 \\
  \delta_3-2 \\
\end{smallmatrix}
\right)$. Clearly, this system of linear equations admits
$\bl x=\frac{1}{3}\left(%
\begin{smallmatrix}
  {\delta_1+\delta_2+\delta_3} \\
  \delta_2+\delta_3-2\delta_1-6 \\
  2\delta_3-\delta_1-\delta_2-6  \\
\end{smallmatrix}%
\right)$ as the only solution. Since $a=2\lambda,
b=d_1^{-1},c=4d_1d_2^{-1}$, it follows that a  necessary
conditions for $\delta_1,\delta_2,\delta_3$ to be eigenvalues of
$\K_{\l,\bl\mu}(0)$ is the inequalities in the statement of the
Lemma.
\end{proof}
\begin{cor}\label{ecur}
Suppose $K^{(\lambda,\bl\mu)}$ and $K^{(\lambda^\prime,\bl\nu)}$
are such that $\K_{\l,\bl\mu}(0)=\K_{\l^\prime,\nu}(0)$ as
matrices, where $(\lambda,\bl\mu)$, $(\lambda^\prime,\bl\nu)$ are
as in Lemma \ref{para}. Then
$(\lambda,\bl\mu)=(\lambda^\prime,\bl\nu)$ as ordered tuples.
\end{cor}
\begin{proof}
Let
$\K_{\l,\bl\mu}(0)=\K_{\l^\prime,\bl\nu}(0)=\mr{diag}(\delta_1,\delta_2,\delta_3)$.
Consider the system of linear equations $A\bl x=\bl b$ and
$A\bl{x^\prime}=\bl{b}$, where $A=\left(%
\begin{smallmatrix}
  1 & -1 & 0 \\
  1 & 1 & -1 \\
  1 & 0 & 1 \\
\end{smallmatrix}%
\right)$, $\bl x=\left(%
\begin{smallmatrix}
  a \\
  b \\
  c \\
\end{smallmatrix}%
\right)$, $\bl{x^\prime}=\left(%
\begin{smallmatrix}
  a^\prime \\
  b^\prime \\
  c^\prime \\
\end{smallmatrix}%
\right)$, $b=\left(%
\begin{smallmatrix}
  \delta_1+2 \\
  \delta_2 \\
  \delta_3-2 \\
\end{smallmatrix}%
\right)$, $a,b,c$ and $a^\prime,b^\prime,c^\prime$ are as in Lemma \ref{para}. As
det$A=3$, $A$ is invertible. Hence $(a,b,c)=(a^\prime,b^\prime,c^\prime)$ and by
Lemma \ref{para}, $(\lambda,\bl\mu)=(\lambda^\prime,\bl\nu)$.
\end{proof}

Suppose $(\delta_1,\delta_2,\delta_3)$, $\delta_i>0$ for
$i=1,2,3$, is given satisfying the inequalities above. Then let us
find $\lambda>1$, $\mu_1,\mu_2>0$ such that $\K_{\l,\bl\mu}(0)=$
diag $(\delta_1,\delta_2,\delta_3)$ with $\bl\mu=(1,\mu_1,\mu_2)$.
We know that  $L(\lambda)\bl\mu^\prime=\bl d$, so $\bl
\mu^\prime=L(\lambda)^{-1}\bl d=\left(%
\begin{smallmatrix}
  1 & 0 & 0 \\
  -\frac{1}{2(\lambda-1)} & 1 & 0 \\
  \frac{1}{\lambda(2\lambda-1)} & -\frac{2}{\lambda} & 1 \\
\end{smallmatrix}%
\right)\left(%
\begin{smallmatrix}
  1 \\
  d_1 \\
  d_2 \\
\end{smallmatrix}%
\right)=\left(%
\begin{smallmatrix}
  1 \\
  d_1-\frac{1}{2(\lambda-1)} \\
  d_2-\frac{2d_1}{\lambda}+\frac{1}{\lambda(2\lambda-1)} \\
\end{smallmatrix}%
\right)$.

Thus
$\mu^2_1=d_1-\frac{1}{2(\lambda-1)}=\frac{1}{b}-\frac{1}{2(\lambda-1)}
= \frac{2(\lambda-1)-b}{2b(\lambda-1)}$ and
$2(\lambda-1)-b=\frac{\delta_1+\delta_2+\delta_3}{3}-2-
\frac{\delta_2+\delta_3-2\delta_1-6}{3}=\delta_1>0$.
$\mu^2_2=d_2-\frac{2d_1}{\lambda}+\frac{1}{\lambda(2\lambda-1)}=\frac{4}{bc}
-\frac{2}{b\lambda}+\frac{1}{\lambda(2\lambda-1)}=\frac{2}{b}(\frac{2}{c}
-\frac{1}{\lambda})+\frac{1}{\lambda(2\lambda-1)}=\frac{2(2\lambda-c)}{bc\lambda}
+\frac{1}{\lambda(2\lambda-1)}=\frac{2(2\lambda-c)(2\lambda-1)+bc}{bc\lambda(2\lambda-1)}
=\frac{2(a-c)(a-1)+bc}{bc\lambda(a-1)}$, where $a,b,c$ are as in
Lemma \ref{curvp}. Thus we have proved the following Theorem.
\begin{thm}\label{ineqthm}
There exists $K^{(\lambda,\bl\mu)}$ such that
$\K_{\l,\bl\mu}(0)=\mr{diag}(\delta_1,\delta_2,\delta_3)$,
$\delta_i>0$ for $i=1,2,3$ if
$$\begin{array}{c}
  \delta_1+\delta_2+\delta_3>6 \\
  \delta_2+\delta_3-2\delta_1>6 \\
  2\delta_3-\delta_1-\delta_2>6  \\
  2(a-c)(a-1)+bc>0 \\
\end{array}$$ where $a,b,c$ are as in Lemma \ref{curvp}.
\end{thm}
\begin{prop} \label{perm1}
Suppose $\delta_i>0$ for $i=1,2,3$ are such that
$\delta_1\neq\delta_2$ and
$2(\delta_1+\delta_2)>\delta_3-6>\mbox{max}\{2\delta_1-\delta_2,2\delta_2-\delta_1\}$.
Then there exists \rks $K^{(\lambda,\bl\mu)}$ and
$K^{(\tilde\lambda,\bl{\tilde\mu})}$ such that
$\K_{\l,\bl\mu}(0)=\mr{diag}(\delta_1,\delta_2,\delta_3)$ and
$\K_{\tilde\l,\tilde{\bl\mu}}(0)=\mr{diag}(\delta_2,\delta_1,\delta_3)$,
where $\lambda,\tilde\lambda>1$, $\bl\mu=(1,\mu_1,\mu_2)$,
$\bl{\tilde\mu}=(1,\tilde\mu_1,\tilde\mu_2)$,
$\mu_\ell,\tilde\mu_\ell>0$ for $\ell=1,2$.
\end{prop}
\begin{proof}
Consider $(\delta_1,\delta_2,\delta_3)$, $\delta_i>0$ for
$i=1,2,3$ such that there exists $K^{(\lambda,\bl\mu)}$ and
$\K_{\l,\bl\mu}(0)=\mr{diag}(\delta_1,\delta_2,\delta_3)$ for some
$\l>1$, $\bl\mu=(1,\mu_1,\mu_2)$ with $\mu_1,\mu_2>0$. So,
$\delta_1,\delta_2,\delta_3$ satisfy the inequalities of Lemma
\ref{nece}. We now produce $\tilde\l>1$,
$\tilde{\bl\mu}=(1,\tilde\mu_1,\tilde\mu_2)$ with
$\tilde\mu_1,\tilde\mu_2>0$ such that
$\K_{\tilde\l,\tilde{\bl\mu}}(0)=\mr{diag}(\delta_2,\delta_1,\delta_3)$.
We recall that $\K_{\tilde\l,\tilde{\bl\mu}}$ is the \C~ of the
metric $\tilde K^{(\tilde\lambda,\tilde{\bl\mu})}(z,z)^t$ and
$\tilde K^{(\tilde\lambda,\tilde{\bl\mu})}$ denotes the
normalization of the \rk $K^{(\tilde\lambda,\tilde{\bl\mu})}$. By
Lemma \ref{curvp} and Remark \ref{tuple} we need to consider the
equations
$$\begin{array}{c}
\tilde a-\tilde b-2=\delta_2\\
\tilde a+\tilde b-\tilde c=\delta_1\\
\tilde a+\tilde c+2=\delta_3 \\
\end{array}$$ where $\tilde a=2\tilde\lambda,\tilde b={\tilde d_1}^{-1},\tilde c=4\tilde d_1{\tilde
d_2}^{-1}$.
This is same as $A\bl {\tilde x}=\bl {\tilde b},$ where $A=\left(%
\begin{smallmatrix}
  1 & -1 & 0 \\
  1 & 1 & -1 \\
  1 & 0 & 1 \\
\end{smallmatrix}%
\right),$ $\bl {\tilde x}=\left(%
\begin{smallmatrix}
  \tilde a \\
  \tilde b \\
  \tilde c \\
\end{smallmatrix}%
\right),$ $\bl {\tilde b} =\left(%
\begin{smallmatrix}
  \delta_2+2 \\
  \delta_1 \\
  \delta_3-2 \\
\end{smallmatrix}%
\right)$. This system of linear equations has  only one
solution, namely,
$\bl x=\frac{1}{3}\left(%
\begin{smallmatrix}
  {\delta_1+\delta_2+\delta_3} \\
  \delta_1+\delta_3-2\delta_2-6 \\
  2\delta_3-\delta_1-\delta_2-6  \\
\end{smallmatrix}%
\right)$. We observe that $a=\tilde a$ and $c=\tilde c$ ~but
$b\neq\tilde b$ if $\delta_1\neq\delta_2$. From Lemma \ref{nece}
and Theorem \ref{ineqthm} we know that there exists
$K^{(\tilde\lambda,\tilde{\bl\mu})}$ such that
$\K_{\tilde\l,\tilde{\bl\mu}}(0)=\mr{diag}(\delta_2,\delta_1,\delta_3)$
if
$$\begin{array}{c}
  \delta_1+\delta_2+\delta_3>6 \\
  \delta_1+\delta_3-2\delta_2>6 \\
  2\delta_3-\delta_1-\delta_2>6  \\
  2(\tilde a-\tilde c)(\tilde a-1)+\tilde b\tilde c>0. \\
\end{array}$$
Hence there exists $K^{(\lambda,\bl\mu)}$ and
$K^{(\tilde\lambda,\tilde{\bl\mu})}$ such that
$\K_{\l,\bl\mu}(0)=$ diag $(\delta_1,\delta_2,\delta_3)$ and
$\K_{\tilde\l,\tilde{\bl\mu}}(0)=\mr{diag}(\delta_2,\delta_1,\delta_3)$
if $$
\begin{array}{c}
 \delta_1+\delta_2+\delta_3>6 \\
  \delta_2+\delta_3-2\delta_1>6 \\
  \delta_1+\delta_3-2\delta_2>6 \\
   2\delta_3-\delta_1-\delta_2>6  \\
   2( a- c)( a-1)+ b c>0 \\
   2(\tilde a-\tilde c)(\tilde a-1)+\tilde b\tilde c>0. \\
\end{array}$$
Suppose $\delta_i>0$ for $i=1,2,3$ are such that
$\delta_1\neq\delta_2$ and
$2(\delta_1+\delta_2)>\delta_3-6>\mbox{max}\{2\delta_1-\delta_2,2\delta_2-\delta_1\}$.
Then we observe that the last inequality implies that
$\delta_2+\delta_3-2\delta_1>6$ and
$\delta_1+\delta_3-2\delta_2>6$, adding these two inequalities we
have $ 2\delta_3-\delta_1-\delta_2>12$. Also
$a-c=\frac{1}{3}(\delta_1+\delta_2+\delta_3)
-\frac{1}{3}(2\delta_3-\delta_1-\delta_2-6)=\frac{1}{3}
(2(\delta_1+\delta_2)-\delta_3)+2$. As $a=\tilde a$ and $c=\tilde
c$, $a-1=\tilde a-1>0$ and the first inequality in the choice of
$\delta_i$ for $i=1,2,3$ implies that $a-c=\tilde a-\tilde c>0$,
so the last two inequalities are satisfied. The first inequality
in the choice of $\delta_i$ for $i=1,2,3$ also implies that
$\delta_3>6$, so the first inequality follows. Hence all the
required inequalities for the existence of $K^{(\lambda,\bl\mu)}$
and $K^{(\tilde\lambda,\bl{\tilde\mu})}$ are satisfied by this
choice of $\delta_i>0$ for $i=1,2,3$.
\end{proof}
\begin{prop} \label{perm2}
Suppose $\delta_i>0$ for $i=1,2,3$ are such that
$\delta_3>\delta_2>3+\frac{\delta_3}{2}$ and
$\delta_1<\mbox{min}\{2\delta_3-\delta_2,2\delta_2-\delta_3\}-6$.
Then there exists \rks $K^{(\lambda,\bl\mu)}$ and
$K^{(\widehat\lambda,\bl{\widehat\mu})}$ such that
$\K_{\l,\bl\mu}(0)=\mr {diag}(\delta_1,\delta_2,\delta_3)$ and
$\K_{\w \l,\w{\bl\mu}}(0)= \mr{diag}(\delta_1,\delta_3,\delta_2)$,
where $\lambda,\w\lambda>1$, $\bl\mu=(1,\mu_1,\mu_2)$,
$\bl{\w\mu}=(1,\w\mu_1,\w\mu_2)$, $\mu_\ell,\w\mu_\ell>0$ for
$\ell=1,2$.
\end{prop}
\begin{proof}
We construct a \rk  $K^{(\w\lambda,\w{\bl\mu})}$ such that $\K_{\w
\l,\w{\bl\mu}}(0)=\mr{diag}(\delta_1,\delta_3,\delta_2)$ for some
$\w\lambda>1$, $\bl{\w\mu}=(1,\w\mu_1,\w\mu_2)$, $\w\mu_\ell>0$
for $\ell=1,2.$ By Lemma \ref{curvp} and Remark \ref{tuple} we
obtain $(\w a,\w b,\w c)$ from  the following set of equations
$$\begin{array}{c}
\w a-\w b-2=\delta_1\\
\w a+\w b-\w c=\delta_3\\
\w a+\w c+2=\delta_2 \\
\end{array}$$ where $\w a=2\w\lambda,\w b={\w d_1}^{-1},\w c=4\w d_1{\w
d_2}^{-1}$.
This is same as $A\bl {\w x}=\bl {\w b},$ where $A=\left(%
\begin{smallmatrix}
  1 & -1 & 0 \\
  1 & 1 & -1 \\
  1 & 0 & 1 \\
\end{smallmatrix}%
\right),$ $\bl {\w x}=\left(%
\begin{smallmatrix}
  \w a \\
  \w b \\
  \w c \\
\end{smallmatrix}%
\right),$ $\bl {\w b} =\left(%
\begin{smallmatrix}
  \delta_2+2 \\
  \delta_1 \\
  \delta_3-2 \\
\end{smallmatrix}%
\right)$. The vector
$\bl x=\frac{1}{3}\left(%
\begin{smallmatrix}
  {\delta_1+\delta_2+\delta_3} \\
  \delta_2+\delta_3-2\delta_1-6 \\
  2\delta_2-\delta_1-\delta_3-6  \\
\end{smallmatrix}%
\right)$ is the only solution of this system  of equations. From
Lemma \ref{nece} and Theorem \ref{ineqthm} we know that there exists
$K^{(\w\lambda,\w{\bl\mu})}$ such that $\K_{\w
\l,\w{\bl\mu}}(0)=\mr{diag}(\delta_1,\delta_3,\delta_2)$  if
$$\begin{array}{c}
  \delta_1+\delta_2+\delta_3>6 \\
  \delta_2+\delta_3-2\delta_1>6 \\
  2\delta_2-\delta_1-\delta_3>6  \\
  2(\w a-\w c)(\w a-1)+\w b\w c>0. \\
\end{array}$$
If $(\delta_1,\delta_2,\delta_3)$, $\delta_i>0$ for $i=1,2,3$ are
such that there exists $K^{(\lambda,\bl\mu)}$ and $\K_{
\l,{\bl\mu}}(0)=\mr{diag}(\delta_1,\delta_2,\delta_3)$. Then
$\delta_i$'s for $i=1,2,3$ satisfies the inequalities of Lemma
\ref{nece}. Hence there exists $K^{(\lambda,\bl\mu)}$ and
$K^{(\w\lambda,\w{\bl\mu})}$ such that
$\K_{\l,\bl\mu}(0)=\mr{diag}(\delta_1,\delta_2,\delta_3)$ and
$\K_{\w \l,\w{\bl\mu}}(0)=\mr{diag}(\delta_1,\delta_3,\delta_2)$
if
$$\begin{array}{c}
   \delta_1+\delta_2+\delta_3>6  \\
  \delta_2+\delta_3-2\delta_1>6 \\
   2\delta_2-\delta_1-\delta_3>6 \\
  2\delta_3-\delta_1-\delta_2>6 \\
  2(a- c)(a-1)+ bc>0 \\
 2(\w a-\w c)(\w a-1)+\w b\w c>0 \\
\end{array}$$
We observe that $a=\w a$ and $b=\w b$ ~but $c\neq\w c$ if
$\delta_2\neq\delta_3$. Suppose $\delta_i>0$ for $i=1,2,3$ are
such that $\delta_3>\delta_2>3+\frac{\delta_3}{2}$ and
$\delta_1<\mbox{min}\{2\delta_3-\delta_2,2\delta_2-\delta_3\}-6$.
Now  the first inequality implies that $\delta_3>6$, hence the
first inequality is satisfied. The last inequality implies that
$2\delta_3-\delta_1-\delta_2>6$ and
$2\delta_2-\delta_1-\delta_3>6$, adding these two inequalities we
have $ \delta_2+\delta_3-2\delta_1>12$. So the first four out of
the set of six inequalities are satisfied. The second, third and
the second, fourth from the set of the six inequalities respectively imply
that $\delta_3-\delta_1>4$ and $\delta_2-\delta_1>4$. An easy
computation involving the expressions for $a,b,c$ and $\w a,\w
b,\w c$ in terms of $\delta_i$ for $i=1,2,3$  shows that $2(a-
c)(a-1)+ bc>0$ and $2(\w a-\w c)(\w a-1)+\w b\w c>0$ are \e~ to
$(\delta_1+\delta_2)(2\delta_1+\delta_2)+\delta_3(\delta_2-\delta_1)+6\delta_1>0$
and
$(\delta_1+\delta_3)(2\delta_1+\delta_3)+\delta_2(\delta_3-\delta_1)+6\delta_1>0$.
These are satisfied as $\delta_2-\delta_1>4$ and
$\delta_3-\delta_1>4$.  Hence all the required inequalities for
the existence of $K^{(\lambda,\bl\mu)}$ and
$K^{(\w\lambda,\bl{\w\mu})}$ are satisfied by this choice of
$\delta_i>0$ for $i=1,2,3$.
\end{proof}
\begin{rem}
The set $\{\delta_i>0 : i=1,2,3\}$ satisfying the  inequalities of
Proposition \ref{perm1} is non-empty.
For instance, take $\delta_1=1$, $\delta_2=2$ and any $\delta_3$ in the interval
$(9,12)$. Then $\{\delta_1,\delta_2,\delta_3\}$ meets the requirement. Similarly,
taking any $\delta_1$ in the interval  $(0,1)$, $\delta_2=7.5$ and
$\delta_3=8$, we find that $\{\delta_1,\delta_2,\delta_3\}$
satisfies the inequalities prescribed in Proposition \ref{perm2}.
Thus the two sets which are obtained from Propositions \ref{perm1} and \ref{perm2}.
are not identical
\end{rem}
\begin{cor}\label{neq}
In Proposition \ref{perm1} and Proposition \ref{perm2},
$(\lambda,\bl\mu)\neq (\tilde\lambda,\tilde{\bl\mu})$ and
$(\lambda,\bl\mu)\neq (\w\lambda,\w{\bl\mu})$.
\end{cor}
\begin{proof}
By Lemma \ref{para} it suffices to show that $(a,b,c)\neq (\tilde
a,\tilde b,\tilde c)$ and $(a,b,c)\neq (\w a,\w b,\w c)$. In
Proposition \ref{perm1} $b\neq \tilde b$ as $\delta_1\neq\delta_2$
and in Proposition \ref{perm2} $c\neq\w c$ as
$\delta_2\neq\delta_3$. 
\end{proof}
Recall that $M^{(\l^\i,\bl\nu)}$ denotes the \mo~ on the \rk~
Hilbert spaces with \rk~ $K^{(\lambda^\i,\bl\nu)}$.
\begin{thm}[\cite{KM}, Theorem 6.2] \label{K}
The reproducing kernels $K^{(\lambda,\bl\mu)}$ and
$K^{(\lambda^\prime,\bl\nu)}$ are \e~ that is, the \mos
$M^{(\l,\bl\mu)}$ and $M^{(\l^\i,\bl\nu)}$  are \u~ if and only if
$(\lambda,\bl\mu)= (\lambda^\i,{\bl\nu})$.
\end{thm}
\begin{cor}\label{only}
Suppose that  $K^{(\lambda,\bl\mu)}$,
$K^{(\tilde\lambda,\tilde{\bl{\mu}})}$ and $K^{(\lambda,\bl\mu)}$,
$K^{(\w\lambda,\w{\bl{\mu}})}$ are as in Proposition \ref{perm1}
and Proposition  \ref{perm2} respectively. Then
\begin{enumerate}
 \item[(a)] the \mos $M^{(\lambda,\bl\mu)}$ and
$M^{(\tilde\lambda,\tilde{\bl{\mu}})}$  are not  \u . \item[(b)]
the \mos $M^{(\lambda,\bl\mu)}$ and $M^{(\w\lambda,\w{\bl\mu})}$
are not \u .
\end{enumerate}
\end{cor}
\begin{proof}
Proof follows immediately from Theorem \ref{K} and Corollary \ref{neq}.
\end{proof}

\begin{rem} \label{perm} In Proposition \ref{perm1} and Proposition
\ref{perm2}, we have shown the following: Given a \rk
$K^{(\lambda,{\bl\mu})}$ such that
$\K_{\l,\bl\mu}(0)=\mr{diag}(\delta_1,\delta_2,\delta_3)$ there
exists a \rk  $K^{(\tilde\lambda,\tilde{\bl\mu})}$ with
$(\lambda,\bl\mu)\neq (\tilde\lambda,\tilde{\bl\mu})$ such that
$\K_{\tilde\l,\tilde{\bl\mu}}(0)=\mr{diag}(\delta_{\rho(1)},\delta_{\rho(2)},\delta_{\rho(3)})$
and given a \rk $K^{(\lambda,{\bl\mu})}$ such that
$\K_{\l,\bl\mu}(0)=\mr{diag}(\delta_1,\delta_2,\delta_3)$ there
exists a \rk $K^{(\w\lambda,\w{\bl\mu})}$ with
$(\lambda,\bl\mu)\neq (\w\lambda,\w{\bl\mu})$ such that $\K_{\w
\l,\w{\bl\mu}}(0)=\mr{diag}(\delta_{\tau(1)},\delta_{\tau(2)},\delta_{\tau(3)})$,
where $\rho,\tau\in \Sigma_3$ with
$\rho(1)=2,\rho(2)=1,\rho(3)=3$, $\tau(1)=1,\tau(2)=3,\tau(3)=2$.
In the next Proposition we prove that there does not exist
$K^{(\lambda^\prime,\bl\nu)}$ with $(\lambda,\bl\mu)\neq
(\lambda^\prime,{\bl\nu})$ such that
$\K_{\l^\prime,\bl\nu}(0)=\mr{diag}(\delta_{\sigma(1)},\delta_{\sigma(2)},\delta_{\sigma(3)})$
if $\sigma\in \Sigma_3$ and $\sigma\neq\rho,\tau$. Obviously,
there exists $K_\sigma^{(\lambda,\bl\mu)}:=P_\sigma
K^{(\lambda,\bl\mu)}P^*_\sigma$ such that $\K_{h_\sigma}(0)=$ diag
$(\delta_{\sigma(1)},\delta_{\sigma(2)},\delta_{\sigma(3)})$ for
all $\sigma\in \Sigma_3$, where $P_\sigma$ is in $\mathcal M(3,\mb C)$  such that $$(P_\sigma)_{ij}=\left\{%
\begin{array}{ll}
    1, & {for~(i,j)=(i,\sigma(i)),} \\
    0, & \hbox{otherwise;} \\
\end{array}%
\right. $$ and $h_\sigma(z)=\tilde
K_\sigma^{(\lambda,\bl\mu)}(z,z)^t$. As the reproducing kernels
$K^{(\lambda,\bl\mu)}$ and $K_\sigma^{(\lambda,\bl\mu)}$ are
equivalent, that is, the multiplication operators on the \rk
Hilbert spaces with reproducing kernels $K^{(\lambda,\bl\mu)}$ and
$K_\sigma^{(\lambda,\bl\mu)}$ are \u,  we do not distinguish
between them. \end{rem}
\begin{prop}\label{exis}
Given a \rk $K^{(\lambda,{\bl\mu})}$ such that
$\K_{\l,\bl\mu}(0)=\mr{diag}(\delta_1,\delta_2,\delta_3)$ there
does not exist \rk $K^{(\lambda^\prime,\bl\nu)}$ with
$(\lambda,\bl\mu)\neq (\lambda^\prime,{\bl\nu})$ such that
$\K_{\l^\prime,\bl\nu}(0)=\mr{diag}(\delta_{\sigma(1)},\delta_{\sigma(2)},\delta_{\sigma(3)})$
if $\sigma\in \Sigma_3$ and $\sigma\neq\rho,\tau$.
\end{prop}
\begin{proof}
Case 1.  Pick $\sigma\in \Sigma_3$ such that
$\sigma(1)=3,\sigma(2)=2,\sigma(3)=1$.

The existence of two \rks  $K^{(\lambda,{\bl\mu})}$ and
$K^{(\lambda^\prime,\bl\nu)}$  such that $\K_{\l,\bl\mu}(0)=\mr{
diag}(\delta_1,\delta_2,\delta_3)$ and
$\K_{\l^\prime,\bl\nu}(0)=\mr{diag}
(\delta_{\sigma(1)},\delta_{\sigma(2)},\delta_{\sigma(3)})$ would
imply, by an application of Lemma \ref{nece} to the ordered
triples $(\delta_1,\delta_2,\delta_3)$ and
$(\delta_{\sigma(1)},\delta_{\sigma(2)},\delta_{\sigma(3)})$,
$$\begin{array}{c}
   \delta_1+\delta_2+\delta_3>6 \\
   \delta_2+\delta_3-2\delta_1>6 \\
  2\delta_3-\delta_1-\delta_2>6   \\
   \delta_{\sigma(1)}+\delta_{\sigma(2)}+\delta_{\sigma(3)}>6 \\
   \delta_{\sigma(2)}+\delta_{\sigma(3)}-2\delta_{\sigma(1)}>6 \\
  2\delta_{\sigma(3)}-\delta_{\sigma(1)}-\delta_{\sigma(2)}>6.  \\
 \end{array}$$ This set of inequalities are \e~ to  $$\begin{array}{c}
  \delta_1+\delta_2+\delta_3>6 \\
   \delta_2+\delta_3-2\delta_1>6 \\
  2\delta_3-\delta_1-\delta_2>6   \\
  \delta_1+\delta_2-2\delta_3>6\\
  2\delta_1-\delta_2-\delta_3>6.\\
\end{array}$$ 
Adding the third and the fourth  from these inequalities
gives $0>12$.

Case 2. Choose $\sigma\in \Sigma_3$  such that
$\sigma(1)=2,\sigma(2)=3,\sigma(3)=1$.

As in the first case the existence of two  \rks
$K^{(\lambda,{\bl\mu})}$ and $K^{(\lambda^\prime,\bl\nu)}$ such
that $\K_{\l,\bl\mu}(0)=\mr{ diag}(\delta_1,\delta_2,\delta_3)$
and $\K_{\l^\prime,\bl\nu}(0)=\mr{diag}
(\delta_{\sigma(1)},\delta_{\sigma(2)},\delta_{\sigma(3)})$ would
imply,  by an application of Lemma \ref{nece} to the ordered
triples $(\delta_1,\delta_2,\delta_3)$ and
$(\delta_{\sigma(1)},\delta_{\sigma(2)},\delta_{\sigma(3)})(=(\delta_2,\delta_3,\delta_1))$,
$$\begin{array}{c}
   \delta_1+\delta_2+\delta_3>6 \\
\delta_2+\delta_3-2\delta_1>6 \\
  2\delta_3-\delta_1-\delta_2>6   \\
  \delta_1+\delta_3-2\delta_2>6 \\
    2\delta_1-\delta_2-\delta_3>6. \\
\end{array}$$ Adding second and fifth of these inequalities gives $0>12$.

Case 3. Take $\sigma\in \Sigma_3$ such that
$\sigma(1)=3,\sigma(2)=1,\sigma(3)=2$.

Finally, continuing in the same manner in the previous two cases,
the existence of two  \rks $K^{(\lambda,{\bl\mu})}$ and
$K^{(\lambda^\prime,\bl\nu)}$ such that $\K_{\l,\bl\mu}(0)=\mr{
diag}(\delta_1,\delta_2,\delta_3)$ and
$\K_{\l^\prime,\bl\nu}(0)=\mr{diag}
(\delta_{\sigma(1)},\delta_{\sigma(2)},\delta_{\sigma(3)})$ would
imply,  by an application of Lemma \ref{nece} to the ordered
triples $(\delta_1,\delta_2,\delta_3)$ and
$(\delta_{\sigma(1)},\delta_{\sigma(2)},\delta_{\sigma(3)})(=(\delta_3,\delta_1,\delta_2))$,
$$\begin{array}{c}
   \delta_1+\delta_2+\delta_3>6 \\
\delta_2+\delta_3-2\delta_1>6 \\
  2\delta_3-\delta_1-\delta_2>6   \\
  \delta_1+\delta_2-2\delta_3>6 \\
    2\delta_2-\delta_3-\delta_1>6. \\
\end{array}$$  Adding third and fourth inequalities from this set of inequalities
we have $0>12$.
\end{proof}
\begin{cor} \label{uniq}
There does not exist any  \mo~   $M^{(\lambda^\i,\bl\nu)}$ not \e~
to $M^{(\lambda,\bl\mu)}$ other than
$M^{(\tilde\l,\tilde{\bl\mu})}$ or $M^{(\w\lambda,\w{\bl\mu})}$
such that $\K_{\l^\i,\bl\nu}(0)=\K_{\l,\bl\mu}(0)$ as sets of
positive numbers, where and $K^{(\lambda,\bl\mu)}$,
$K^{(\tilde\lambda,\tilde{\bl\mu})}$, $K^{(\w\lambda,\w{\bl\mu})}$
are as in Proposition \ref{perm1} and Proposition \ref{perm2}.
\end{cor}
\begin{proof}
Combining Corollary \ref{only}, Theorem \ref{K}, Corollary
\ref{ecur} and Proposition \ref{exis}, we obtain a proof of this
corollary.
\end{proof}
\begin{rem} \label{r2}
We discuss the case $m=1$. Proceeding as in Lemma \ref{curvp} we
see that $\K_{\l,\bl\mu}(0)=\mr{diag}(a-b-1,a+b+1)$, where
$\lambda>1/2$, $\bl\mu=(1,\mu_1)$, $\mu_1>0$, $a=2\lambda$,
$b={d_1}^{-1}$, $d_1$ is defined as before. If
$\K_{\l,\bl\mu}(0)=\mr{diag}(\delta_1,\delta_2)$, $\delta_i>0$ for
$i=1,2$, for some $\l>1/2$ and $\bl\mu=(1,\mu_1)$, $\mu_1>0$. Then
arguing as in Lemma \ref{nece} one notes that
$a=2\l=\frac{\delta_1+\delta_2}{2}$,
$b=\frac{\delta_2-\delta_1-2}{2}$. As $a=2\l>1$ and $b=d^{-1}>0$
it follows that $\delta_1+\delta_2>2$ and $\delta_2-\delta_1>2$
are necessary conditions for existence of a \rk
$K^{(\l,\bl\mu)}$such that
$\K_{\l,\bl\mu}(0)=\mr{diag}(\delta_1,\delta_2)$. If $\delta_i>0$
for $i=1,2$, proceeding as in Theorem \ref{ineqthm}, one observes
that $\delta_2-\delta_1>2$, $\delta_1+\delta_2>2$ and
$d_1>\frac{1}{2\lambda-1}=\frac{2}{\delta_1+\delta_2-2}$ are the
sufficient conditions for existence of a  \rk~
$K^{(\lambda,\bl\mu)}$ such that  $\K_{\l,\bl\mu}(0)=\mr{diag}
(\delta_1,\delta_2)$. Conversely, if $\delta_i>0$ for $i=1,2$ and
$\delta_2-\delta_1>2$ then clearly $\delta_1+\delta_2>2$ and
$d_1=\frac{2}{\delta_2-\delta_1-2}>\frac{2}{\delta_1+\delta_2-2}$.
So $\delta_i>0$ for $i=1,2$ and $\delta_2-\delta_1>2$ are the
necessary and sufficient conditions for the existence of  \rk~
$K^{(\lambda,\bl\mu)}$ such that
$\K_{\l,\bl\mu}(0)=\mr{diag}(\delta_1,\delta_2)$.
\end{rem}
\begin{rem}
If $\delta_i>0$ for $i=1,2$ and $\delta_2-\delta_1>2$ there does
not exist a  \rk $K^{(\lambda^\prime,\bl\nu)}$ other than
$K^{(\lambda,\bl\mu)}$ (up to \eq~ as discussed in Remark
\ref{perm}) such that $\K_{\l^\prime,\bl\nu}(0)=
\mr{diag}(\delta_2,\delta_1)$. If $K^{(\lambda^\prime,\nu)}$
exists satisfying the above requirements then from Remark \ref{r2}
we see that both of $\delta_2-\delta_1>2$ and
$\delta_1-\delta_2>2$ have to be simultaneously satisfied. This is
impossible as they imply $0>4$. Hence there does not exist in\e~
\mos $M^{(\lambda,\bl\mu)}$ and $M^{(\lambda^\prime,\bl\nu)}$ such
that $\K_{\l,\bl\mu}(0)=\K_{\l^\i,\bl\nu}(0)$ as sets.
\end{rem}
\begin{thm} \label{ineqi}
Suppose that $K^{(\l,\bl\mu)}$, $K^{(\tilde\l,\tilde{\bl\mu})}$
and $K^{(\l,\bl\mu)}$, $K^{(\w\l,\w{\bl\mu})}$ are as in
Proposition \ref{perm1} and Proposition \ref{perm2} respectively.
Then
\begin{enumerate}
\item[(i)] the \mos $M^{(\l,\bl\mu)}$ and
$M^{(\tilde\l,\tilde{\bl\mu})}$  are not \e~ although
$\K_{\l,\bl\mu}(z)$ and $\K_{\t\l,\t{\bl\mu}}(z)$ are \u~ for $z$
in $\D$. \item[(ii)] the \mos $M^{(\l,\bl\mu)}$ and
$M^{(\tilde\l,\tilde{\bl\mu})}$  are not \e~ although
$\K_{\l,\bl\mu}(z)$ and $\K_{\w\l,\w{\bl\mu}}(z)$ are \u~ for $z$
in $\D$.
\end{enumerate}
\end{thm}
\begin{proof}
From \ref{perm1} we see that the \Cs of the associated bundles
have the same set of eigenvalues at zero namely,
$\{\delta_1,\delta_2,\delta_3\}$. Since curvature is self-adjoint
the set of eigenvalues is the complete set of unitary invariants
for the curvature. So, $\K_{\l,\bl\mu}(0)$ and $\K_{\tilde
\l,\tilde{\bl\mu}}(0)$ are unitarily equivalent. As the operators
$M^{(\l,\bl\mu)}$ and  $M^{(\tilde\l,\tilde{\bl\mu})}$ are \H, by
an application of Theorem \ref{hom} we see that
$\K_{\l,\bl\mu}(z)$ and $\K_{\tilde \l,\t{\bl\mu}}(z)$ are \ue~
for $z\in\D$. Now (i) follows from part (a) of Corollary
\ref{only}. The proof of (ii) of this theorem is similar.
\end{proof}
The proof of the next Theorem will be completed after proving a
sequence of Lemmas. We omit the easy proof of the first of these lemmas.
\begin{thm} \label{com}
Suppose that $M^{(\l,\bl\mu)}$ and $M^{(\l^\i,\bl\nu)}$ are not
\u~ and the two \Cs  $\K_{\l,\bl\mu}(z)$ and
$\K_{\l^\i,\bl\nu}(z)$ are \u~ for $z$ in $\D$. Then there does
not exist any invertible matrix $L$ in $\mathcal M(3,\mb C)$
satisfying $L\K_{\l,\bl\mu}(0)=\K_{\l^\i,\bl\nu}(0)L$ for which
$L(\K_{\l,\bl\mu})_{\bar z}(0)=(\K_{\l^\i,\bl\nu})_{\bar z}(0)L$
also. In other words, the covariant derivative of order $(0,1)$
detects the inequivalence.
\end{thm}
\begin{lem}\label{mat}
Suppose that
$\Delta=\big(\!\!\big(k_i\delta_{ij}\big)\!\!\big)_{i,j=1}^n$,
$\Delta_\sigma=\big(\!\!\big(k_{\sigma(i)}\delta_{ij}\big)\!\!\big)_{i,j=1}^n$,
$k_i\neq k_j$~ if $~i\neq j$ and
$C=\big(\!\!\big(c_{ij}\big)\!\!\big)_{i,j=1}^n$ in $\mathcal
M(n,\mb C)$ is such that $C\Delta=\Delta_\sigma C$. Then
$c_{ij}=c_{ij}\delta_{\sigma(i),j}$ for $i,j=1,\ldots,n$, where
$\sigma$ is in $S_n$, $S_n$ denotes the permutation group of
degree $n$.
\end{lem}
\begin{cor}  \label{mat1}
Suppose that there exists \rks $K^{(\lambda,\bl\mu)}$,
$K^{(\tilde\lambda,\bl{\tilde\mu})}$ such that
$\K_{\l,\bl\mu}(0)=\mr{diag}(\delta_1,\delta_2,\delta_3)$,
$\K_{\tilde\l,\tilde{\bl\mu}}(0)=\mr{diag}(\delta_{\rho(1)}
,\delta_{\rho(2)},\delta_{\rho(3)})$, $\delta_1\neq\delta_2$ and
$C=\big(\!\!\big(c_{ij}\big)\!\!\big)_{i,j=1}^3$ in $\mathcal
M(3,\mb C)$ is such that
$C\K_{\l,\bl\mu}(0)=\K_{\tilde\l,\tilde{\bl\mu}}(0)C$, then
$c_{ij}=c_{ij}\delta_{\rho(i),j}$ for $i,j=1,2,3$, where $\rho\in
\Sigma_3$ is such that $\rho(1)=2,\rho(2)=1,\rho(3)=3$.
\end{cor}
\begin{proof}
The proof of this Corollary is immediate from Lemma \ref{mat} once
we ensure that $\delta_1,\delta_2,\delta_3$ are distinct.
Recalling notations from Lemma \ref{curvp} we write
$\delta_1=a-b-2$, $\delta_2=a+b-c$, $\delta_3=a+c+2$. Clearly,
$\delta_3-\delta_1=b+c+4>0$.  Recalling notations from Proposition
\ref{perm1} one has $\delta_2=\t a-\t b-2$, $\delta_1=\t a+\t b-\t
c$, $\delta_3=\t a+\t c+2$. So, $\delta_3-\delta_2=\t b+\t c+4>0$.
We have $\delta_3>\delta_1$, $\delta_3>\delta_2$ and $\delta_1\neq
\delta_2$ by hypothesis. Hence the proof is complete.
\end{proof}
The proof of the next Corollary is similar.
\begin{cor} \label{mat2}
Suppose that there exists \rks $K^{(\lambda,\bl\mu)}$,
$K^{(\w\lambda,\bl{\w\mu})}$ such that
$\K_{\l,\bl\mu}(0)=\mr{diag}(\delta_1,\delta_2,\delta_3)$,
$\K_{\w\l,\w{\bl\mu}}(0)=\mr{diag}(\delta_{\tau(1)}
,\delta_{\tau(2)},\delta_{\tau(3)})$, $\delta_2\neq\delta_3$ and
$C=\big(\!\!\big(c_{ij}\big)\!\!\big)_{i,j=1}^3$ in $\mathcal
M(3,\mb C)$ is such that
$C\K_{\l,\bl\mu}(0)=\K_{\w\l,\w{\bl\mu}}(0)C$, then
$c_{ij}=c_{ij}\delta_{\tau(i),j}$ for $i,j=1,2,3$, where $\tau\in
\Sigma_3$ is such that $\tau(1)=1,\tau(2)=3,\tau(3)=2$.
\end{cor}

\begin{lem} \label{det}
Suppose that
$C=\big(\!\!\big(c_{ij}\delta_{\sigma(i),j}\big)\!\!\big)_{i,j=1}^3$
 for $\sigma=\rho, \tau$ in $\Sigma_3$. Then $C$ is invertible if and
 only if $c_{i,\sigma(i)}\neq 0$ for $i=1,2,3$ and $\sigma=\rho, \tau$ in
$\Sigma_3$.
\end{lem}
\begin{proof}
Since det $C\neq 0$ if and only if $c_{i,\sigma(i)}\neq 0$ for
$i=1,2,3$ and $\sigma=\rho, \tau$ in $\Sigma_3$, the proof is
complete.
\end{proof}
The proof of the following Lemma is straight forward. We recall
that $\big(\!\!\big(S_m(c_1,\ldots,c_m)\big)\!\!\big)_{\ell,
p}=c_\ell\delta_{p+1,\ell}$, $0\leq p,\ell\leq m$.
\begin{lem} \label{det1}
Suppose that
$C=\big(\!\!\big(c_{ij}\delta_{\sigma(i),j}\big)\!\!\big)_{i,j=1}^3$,
$c_{i,\sigma(i)}\neq 0$ for $i=1,2,3$ and $\sigma=\rho, \tau$ in
$\Sigma_3$ is such that $CS_2(c_1,c_2)^t=S_2(\t c_1,\t c_2)^tC$
for $c_i,\t c_i$ in $\mb C$ for $i=1,2$. Then $c_i=\t c_i=0$ for
$i=1,2$.
\end{lem}
\begin{lem} \label{non}
$(\K_{\l,\bl\mu})_{\bar z}(0)$ is not the zero matrix.
\end{lem}
\begin{proof}
If possible let $(\K_{\l,\bl\mu})_{\bar z}(0)=0$. Then it follows
from Lemma \ref{curvp} that  $-\sqrt b(1+b-\frac{c}{2})=-\sqrt
c(1+c-\frac{b}{2})=0$. Equivalently,
$1+b-\frac{c}{2}=1+c-\frac{b}{2}$, as $b$ and $c$ are positive.
This implies that $b=c$. So, $(\K_{\l,\bl\mu})_{\bar z}(0)=0$
implies by an application of Lemma \ref{curvp} that $-\sqrt
b(1+\frac{b}{2})=0$, which is impossible as $b$ is positive.
\end{proof}
\begin{proof}[Proof of Theorem \ref{com}:] We observe by applying Proposition
\ref{perm1}, Proposition \ref{perm2} and Proposition \ref{exis}
that if $M^{(\l^\i,\bl\nu)}$ is an \mo not \u~ to
$M^{(\l,\bl\mu)}$ then $(\l^\i,\bl\nu)=(\t\l,\t{\bl\mu})$ or
$(\w\l,\w{\bl\mu})$. We reach the conclude by an straight forward
application of Corollary \ref{mat1}, Corollary \ref{mat2}, Lemma
\ref{det}, Lemma \ref{det1} and Lemma \ref{non}.
\end{proof}
The comments below are for the class of \H~ operators constructed
in \cite{KM}.
\begin{rem}
Unfortunately, although we are able to carry out similar
calculations for these operators of rank $\geq 4$, it is not
clear, if this would completely answer the question raised in
\cite[page. 39]{CD2}. Indeed, for rank $3$, we have shown that the
simultaneous \ue~ class of the \C~ at $0$ along with  the
covariant derivative of \C~ at $0$ of order $(0,1)$ determines the
\ue~ class of these operators. Similarly,  for rank $2$, the \ue~
class of the \C~ at $0$  determines the \ue~ class of the
operator.
\end{rem}

\end{document}